\documentclass[amsthm]{elsart}

\usepackage{yjsco}
\usepackage{natbib}
\usepackage{amsthm, amsmath}   % For Latex2e
\usepackage{latexsym,amssymb}
\usepackage{savesym}
\savesymbol{AND}
\usepackage{algorithmic}
\usepackage{enumitem}

\usepackage{amssymb}
\usepackage[usenames,dvipsnames]{pstricks}
\usepackage{epsfig}
\usepackage{pst-node}
\usepackage{pst-grad} % For gradients
\usepackage{pst-plot} % For axes
\usepackage{pst-3dplot}

\def\QED{\hfill$\Box$}

\def\P{{\mathbb P}}
\def\N{{\mathbb N}}
\def\A{{\mathbb A}}
\def\Z{{\mathbb Z}}
\def\Q{{\mathbb Q}}
\def\R{{\mathbb R}}

\def\gcd{{\rm{gcd}}}

\def\IA{I_{\mathcal A}}
\def\IB{I_{\mathcal B}}
\def\mA{\mathcal A}
\def\mB{\mathcal B}
\def\mC{\mathcal C}
\def\mH{\mathcal H}
\def\x{\mathbf{x}}
\def\t{\mathbf{t}}
\def\kx{k[\x]}
\def\kt{k[\t]}

\newtheorem{Theorem}{Theorem}[section]
\newtheorem{Lemma}[Theorem]{Lemma}
\newtheorem{Corollary}[Theorem]{Corollary}
\newtheorem{Proposition}[Theorem]{Proposition}
\newtheorem{Remark}[Theorem]{Remark}
\newtheorem{Example}[Theorem]{Example}

\newtheorem{Definition}[Theorem]{Definition}

\def\QED{\hfill$\Box$}

\setitemize{leftmargin=.9cm} \setlist[1]{itemsep=6pt}

\begin{document}
\parskip 8pt

\begin{frontmatter}

\title{Complete intersections in simplicial toric varieties }

\thanks{Both authors were partially supported by Ministerio de
Ciencia e Innovaci\'on,  Spain (MTM2010-20279-C02-02). }

\author{Isabel Bermejo}
\address{Facultad de Matem\'aticas, Universidad de La Laguna, 38200 La Laguna, Tenerife,
Spain} \ead{ibermejo@ull.es}

\author{Ignacio Garc\'{\i}a-Marco}\footnote{Corresponding author. Tel: (+34) 922 31 81 46}
\address{Facultad de Matem\'aticas, Universidad de La Laguna, 38200 La Laguna, Tenerife,
Spain} \ead{iggarcia@ull.es}

\begin{abstract}
Given a set $\mA = \{a_1,\ldots,a_n\} \subset \N^m$ of nonzero
vectors defining a simplicial toric ideal $\IA \subset
k[x_1,\ldots,x_n]$, where $k$ is an arbitrary field, we provide an
algorithm for checking whether $\IA$ is a complete intersection.
This algorithm does not require the explicit computation of a
minimal set of generators of $\IA$. The algorithm is based on the
application of some new results concerning toric ideals to the
simplicial case. For homogenous simplicial toric ideals, we provide
a simpler version of this algorithm. Moreover, when $k$ is an
algebraically closed field, we list all ideal-theoretic complete
intersection simplicial projective toric varieties that are either
smooth or have one singular point.
\end{abstract}


\end{frontmatter}

\section{Introduction}\label{introduction}
Let $k$ be an arbitrary field and $\kx = k[x_1,\ldots,x_n]$ and $\kt
= k[t_1,\ldots,t_m]$ two polynomial rings over $k$. A {\it binomial}
in a polynomial ring is a difference of two monomials. Let $\mA
=\{a_1,\ldots,a_n\}$ be a set of nonzero vectors in $\N^m$; each
vector $a_i = (a_{i1},\ldots,a_{im})$ corresponds to a monomial
$\t^{a_i} = t_1^{a_{i1}} \cdots t_m^{a_{im}}$ in $\kt$. The {\it
toric set} $\Gamma$ determined by $\mA$ is the subset of the affine
space $\A_k^n$ given parametrically by $x_i = u_1^{a_{i1}}\cdots
u_m^{a_{im}}$ for all $i \in \{1,\ldots,n\}$, i.e.,
$$
\Gamma=\{(u_1^{a_{11}}\cdots u_m^{a_{1m}},\ldots,u_1^{a_{n1}}\cdots
u_m^{a_{nm}})\in\mathbb{A}_k^n\, \vert\, u_1,\ldots,u_m \in k \}\,.
$$

The kernel of the homomorphism of $k$-algebras $\varphi\colon
\kx\rightarrow \kt;\ x_i\longmapsto \t^{a_i}$ is called the {\it
toric ideal\/} of $\Gamma$ and will be denoted by
 $\IA$. For every $b = (b_1,\ldots,b_n) \in
\N^n$ one sets the $\mA$-degree of the monomial $\x^b  \in \kx$ as
${\rm deg}_{\mA}(\x^b) := b_1 a_1 + \cdots + b_n a_n \in \N^m$. One
says that a polynomial $f \in \kx$ is $\mA$-homogeneous if its
monomials have the same $\mA$-degree. By
\cite[Corollary~4.3]{Sturm}, it is an $\mA$-{\it homogeneous
binomial\/} ideal, i.e., $\IA$ is generated by $\mA$-homogeneous
binomials. According to \cite[Lemma 4.2]{Sturm}, the height of $\IA$
is equal to $n-{\rm dim}(\Q \mA)$. By
\cite[Corollary~7.1.12]{monalg}, if $k$ is an infinite field, $\IA$
is the ideal $I(\Gamma)$ of the polynomials vanishing in $\Gamma$.
The variety $V(\IA) \subset \A_k^n$ is called an {\it affine toric
variety}.

The ideal $\IA$ is a {\it complete intersection} if $\mu(\IA) = {\rm
ht}(\IA)$, where $\mu(\IA)$ denotes the minimal number of generators
of $\IA$. Equivalently, $\IA$ is a complete intersection if there
exists a set of $s = n-{\rm dim}(\Q \mA)$ $\mA$-homogeneous
binomials $g_1,\ldots,g_{s}$ such that $\IA=(g_1,\ldots,g_{s}).$ The
problem of determining complete intersection toric ideals has a long
history; see the introduction of \citet{MT} and the references
there.

We denote by ${\rm Cone}(\mA)$ the cone spanned by $\mA$, i.e.,
${\rm Cone}(\mA) = \{ \sum_{i = 1}^n \alpha_i a_i \, \vert \,
\alpha_i \in \R_{\geq 0} \}$. An {\it extreme ray} of ${\rm
Cone}(\mA)$ is a set $F_a := {\rm Cone}(\{a\})$ such that whenever
$x,y \in {\rm Cone}(\mA)$ satisfy that $x + y \in F_a$, then $x,y
\in F_a$. It is well known, see for example \cite[Lemma
1.2.15]{CLS}, that a set $\{a_1',\ldots,a_s'\} \subset \R^m$ is a
minimal set of generators of ${\rm Cone}(\mA)$ if and only if the
extremal rays of ${\rm Cone}(\mA)$ are $F_{a_1'},\ldots,F_{a_s'}$.
Thus, the number of extremal rays of ${\rm Cone}(\mA)$ is $\geq {\rm
dim}(\Q \mA)$. When equality holds the toric ideal $\IA$ is said to
be a {\it simplicial toric ideal} and $V(\IA)$ an {\it affine
simplicial toric variety}. If $\IA$ is homogeneous, then $V(\IA)
\subset \P_k^{n-1}$ is called a {\it simplicial projective toric
variety}.

The aim of this work is to obtain and implement an efficient
algorithm for checking whether a simplicial toric ideal is a
complete intersection that does not require the explicit computation
of a minimal set of generators of the ideal.

This work follows the line we began in \cite{BGS}, where we proposed
an algorithm for checking whether the toric ideal of an affine
monomial curve is a complete intersection. That algorithm was based
on the ideas introduced in \cite{BGRV} and we implemented it in {\sc
Singular} \citep{DGPS}, giving rise to the library {\tt cimonom.lib}
\citep{BGSlib}. This work is a non trivial generalization of
\cite{BGS} for simplicial toric ideals and gives rise to the {\sc
Singular} library {\tt cisimplicial.lib} \citep{BGlib2}, which
generalizes, outperforms and substitutes our previous {\tt
cimonom.lib}.

For proving correctness of our algorithm we will use the following
direct consequence of \cite[Theorem 1.4]{Rosalesgluing}.
\begin{quote} If $\mA$ is a gluing of $\mA_1$ and $\mA_2$ and both
$I_{\mA_1}$ and $I_{\mA_2}$ are complete intersections, then so is
$\IA$. \end{quote} Recall that $\mA$ is a gluing of $\mA_1$ and
$\mA_2$ if $\mA = \mA_1 \sqcup \mA_2$ and there exists $\alpha \in
\N \mA_1 \cap \N \mA_2$ such that $\Z \alpha = \Z \mA_1 \cap \Z
\mA_2$.

It is worth pointing out that Fischer, Morris and Shapiro provided
in \citep{FMS} a theoretical characterization of the property of
being a complete intersection in toric ideals by proving that
whenever $\IA$ is a complete intersection, there exist $\mA_1, \mA_2
\subset \mA$ such that $\mA$ is a gluing of $\mA_1$ and $\mA_2$ and
both $I_{\mA_1}, I_{\mA_2}$ are complete intersections. This result
was improved in \cite{GR} for the particular case of simplicial
toric ideals. Our approach to the problem of determining complete
intersection simplicial toric ideals is different in nature to that
of \cite{FMS} and \cite{GR}.

The main achievement of this work is Algorithm CI-simplicial,  which
receives as input any set $\mA \subset \N^m$ such that $\IA$ is a
simplicial toric ideal and returns {\sc True} if $\IA$ is a complete
intersection or {\sc False} otherwise. Moreover, whenever $\IA$ is a
complete intersection, the algorithm provides without any extra
effort a minimal set of $\mA$-homogeneous generators of $\IA$. This
algorithm is based on the application of some new results concerning
complete intersection toric ideals, namely Theorems
\ref{teoremaReduccion} and \ref{ci-iguales}, to the simplicial case.
Correctness of Algorithm CI-simplicial is proved in Theorem
\ref{algoritmosimplicial}, which is the main result of this paper.

The structure of the paper is the following.

Sections \ref{sec2} and \ref{sec3} are devoted to present two
techniques on complete intersection toric ideals. In section
\ref{sec2} we prove Theorem \ref{teoremaReduccion}. This result is
based on the idea of associating to $\mA$ another set $\mA_{red}
\subset \N^m$ that can be either empty or defines a toric ideal
$I_{\mA_{red}}$ satisfying that $\IA$ is a complete intersection if
and only if either $\mA_{red} = \emptyset$ or $I_{\mA_{red}}$ is a
complete intersection. Moreover, in case $\mA_{red}$ is not empty,
$I_{\mA_{red}} \subset k[x_1,\ldots,x_t]$ with $t \leq n$ and the
degrees of the generators of $I_{\mA_{red}}$ are lower than those in
$\IA$. The construction of $\mA_{red}$ is described algorithmically
and, as a consequence of this, if $\mA_{red} = \emptyset$ or one
knows a minimal set of generators of $I_{\mA_{red}}$, one can
recover a minimal set of generators of $\IA$.

In Section \ref{sec3} we prove Theorem \ref{ci-iguales}. For some $i
\in \{1,\ldots,n\}$, one can set
\begin{center}$m_i := {\rm min}\left\{b \in \Z^+ \, \vert \, b a_i
\in \sum_{j \in \{1,\ldots,n\} \atop j \neq i} \N a_j\right\}$
\end{center} and in case there exist $i,j:\, 1 \leq i < j \leq n$ such
that $m_i a_i = m_j a_j$, we define $a_i^{\,\prime} := (a_{i
1}^{\,\prime},\ldots,a_{i m}^{\,\prime}) \in \N^m$, where $a_{i
k}^{\,\prime} = 0$ if $a_{i k} = 0$ or $a_{i k}^{\,\prime} =
\gcd\{a_{i k}, a_{j k}\}$ otherwise. Theorem \ref{ci-iguales} states
among other things that whenever $m_i a_i = m_j a_j$ for some $i,j:
\, 1 \leq i < j \leq n$, if $\IA$ is a complete intersection, then
so is $I_{\mA_{(i,j)}}$, where $\mA_{(i,j)} := (\mA \setminus
\{a_i,a_j\}) \cup \{a_i^{\,\prime}\}$.  Moreover, the toric ideal
$I_{\mA_{(i,j)}}$ belongs to a ring of polynomials with $n-1$
variables and its height is one unit less than $\IA$. For proving
this result we use the relationship between mixed dominating
matrices and complete intersection toric ideals established in
\cite{Fisher-Shapiro}.

Section \ref{sec4} is devoted to design Algorithm CI-simplicial, the
main result of this section is Theorem \ref{algoritmosimplicial},
where correctness of the algorithm is proved. A key result for
obtaining the algorithm is Proposition \ref{simplicial}, which
asserts that if $\IA$ is a complete intersection simplicial toric
ideal, then either there exist $i,j$ such that $m_i a_i = m_j a_j$
or $\mA_{red} = \emptyset$. Moreover, if $\IA$ is a complete
intersection and $m_i a_i = m_j a_j$, then by Theorem
\ref{ci-iguales} the set $\mA_{(i,j)}$ determines another complete
intersection simplicial toric ideal $I_{\mA_{(i,j)}}$. The idea
under Algorithm CI-simplicial is to apply Theorem \ref{ci-iguales}
as many times as possible, until we get a set $\mB$ such that $\IB$
is a simplicial toric ideal and $\IB$ is a complete intersection if
$\IA$ so is. Then we compute $\mB_{red}$. If $\mB_{red} \not=
\emptyset$, then $\IA$ is not a complete intersection. Nevertheless,
if $\mB_{red} = \emptyset$ some extra conditions have to be verified
in order to determine whether $\IA$ is a complete intersection.
These conditions consist of checking whether some elements belong to
certain subsemigroups of $\N^m$.

In Section \ref{sec5} we study the complete intersection property
for homogeneous simplicial toric ideals. Firstly,  as a direct
consequence of Theorem \ref{teoremaReduccion} and Proposition
\ref{simplicial}, we get in Corollary \ref{coralg} that a
homogeneous simplicial toric ideal $\IA$ is a complete intersection
if and only if  $\mA_{red} = \emptyset$. This result provides a
simpler version of Algorithm CI-simplicial for homogeneous
simplicial toric ideals that only consists of computing $\mA_{red}$
and checking if $\mA_{red} = \emptyset$. As non-trivial consequences
of Corollary \ref{coralg}, when $k$ is an algebraically closed
field, we prove in Theorem \ref{smoothIC} that there is only one
smooth simplicial projective toric variety that is an
ideal-theoretic complete intersection, which is the projective
monomial curve defined parametrically by $x_1 = u_1^2,\, x_2 =
u_2^2,\, x_3 = u_1 u_2$. Moreover, in  Theorem \ref{1singIC}, we
list all simplicial projective toric varieties having one singular
point that are ideal-theoretic complete intersection. Recall that a
variety is an ideal-theoretic complete intersection if its defining
ideal is a complete intersection.

Finally, in Section \ref{sec6}, we describe the implementation
details of the algorithms for determining whether a homogeneous
simplicial toric ideal and a simplicial toric ideal is a complete
intersection. We have implemented these algorithms in C++ and in
{\sc Singular}. Our implementation in {\sc Singular} gave rise to
the distributed library {\tt cisimplicial.lib} \citep{BGlib2}.
Computational experiments show that our implementation is able to
solve large-size instances.

\section{From $\IA$ to $I_{\mA_{red}}$}\label{sec2}

Starting from $\mA \subset \N^m$, we are going to construct another
set $\mA_{red} \subset \N^m$ such that $\IA$ is a complete
intersection if and only if $\mA_{red} = \emptyset$ or
$I_{\mA_{red}}$ is a complete intersection. Moreover, if $\mA_{red}
\neq \emptyset$, then $\mA_{red} = \{a_1',\ldots,a_t'\}$ with $t
\leq n$ and the generators of $I_{\mA_{red}}$ have lower degrees
than those of $\IA$.

In order to explain how to construct $\mA_{red}$ from $\mA$, we need
Lemmas \ref{redpertenece}, \ref{red2}, which are easy to prove, and
Proposition \ref{red}.

\begin{Lemma}\label{redpertenece}If there exist $\{\alpha_j\}_{j \in \{1,\ldots,n\} \atop j \neq i} \subset \N$
such that $a_i = \sum_{j \in \{1,\ldots,n\} \atop j \neq i} \alpha_j
a_j$ for some $i \in \{1,\ldots,n\}$, then
\begin{itemize}  \item[{\rm (a)}] $\IA = I_{\mA \setminus \{a_i\}}
\cdot \kx + (x_i - \prod_{j \in \{1,\ldots,n\} \atop j \neq i}
x_j^{\alpha_j})$.
\item[{\rm (b)}]
$\mu(\IA) = \mu(I_{\mA \setminus \{a_i\}}) + 1$.
\item[{\rm (c)}] $\IA$ is a complete intersection $\Longleftrightarrow$ $I_{\mA
\setminus \{a_i\}}$ is a complete intersection. \end{itemize}
\end{Lemma}

\smallskip

\begin{Lemma}\label{red2}If $a_i \notin \sum_{j \in \{1,\ldots,n\} \atop j \neq i}
\Q\, a_j$ for some $i \in \{1,\ldots,n\}$, then
\begin{itemize} \item[{\rm (a)}] $\IA = I_{\mA \setminus \{a_i\}}
\cdot \kx$.
\item[{\rm (b)}]
$\mu(\IA) = \mu(I_{\mA \setminus \{a_i\}})$.
\item[{\rm (c)}] $\IA$ is a complete intersection $\Longleftrightarrow$ $I_{\mA
\setminus \{a_i\}}$ is a complete intersection. \end{itemize}
\end{Lemma}

Take $i \in \{1,\ldots,n\}$ such that $a_i \in  \sum_{j \in
\{1,\ldots,n\}\atop j \neq i} \Q\, a_j$, we define
\begin{center} $B_i := {\rm min}\{b \in \Z^+\, \vert \ b\, a_i \in
 \sum_{j \in \{1,\ldots,n\} \atop j \neq i} \Z\, a_j \}$
\end{center}
and we get the following result, which is a generalization of the
result for $1$-dimensional toric ideals \cite[Lemma 1.2]{Mor2}, see
also \cite[Lemma 3.2]{Mor}. Since the idea under its proof is
analogue to that of \cite[Lemma 1.2]{Mor2} we do not include it
here.

\begin{Proposition}\label{red}Set $\mA^{\,\prime} :=
\{a_1,\ldots,a_{i-1}, b\,a_i,a_{i+1}, \ldots,a_n\} \subset \N^m$,
where $b \in \Z^+$ is a divisor of $B_i$. Then,
\begin{itemize} \item[{\rm (a)}] $\IA
= \rho(I_{\mA^{\,\prime}}) \cdot \kx$, where $\rho: \kx
\longrightarrow \kx$ is the $k$-homomorphism defined by $x_i \mapsto
x_i^{b}$\,, $x_j \mapsto x_j$\, for $j \in \{1,\ldots,n\} \setminus
\{i\}$.
 \item[{\rm (b)}]  $\mu(\IA) = \mu(I_{\mA'})$.
\item[{\rm (c)}] $\IA$ is a complete intersection
$\Longleftrightarrow$ $I_{\mA'}$ is a complete intersection.
\end{itemize}
\end{Proposition}

Applying Lemmas \ref{redpertenece} and \ref{red2} and Proposition
\ref{red} as many times as possible, we associate to $\mA$ a subset
$\mA_{red} \subset \N^m$ which can be either empty or $\mA_{red} =
\{a_1^{\,\prime},\ldots,a_t^{\,\prime}\}$ and for every $i \in
\{1,\ldots,t\}$ the following hold:
\begin{itemize}
\item $ a_i^{\,\prime} \notin \sum_{j \in \{1,\ldots,t\}\atop j \neq i}
\N\, a_j^{\,\prime}$,
\item $ a_i^{\,\prime} \in \sum_{j \in \{1,\ldots,t\}\atop j \neq i}
\Q\, a_j^{\,\prime}$, and \item if $B_i^{\, \prime} := {\rm min}\{b
\in \Z^+\,\vert\, b a_i^{\,\prime} \in \sum_{j \in \{1,\ldots,t\}
\atop j \neq i} \Z a_j^{\,\prime} \}$, then $B_i^{\,\prime} = 1$
($a_i^{\,\prime} \in \sum_{j \in \{1,\ldots,t\}\atop j \neq i} \Z\,
a_j^{\,\prime}).$
\end{itemize}

\medskip

As a direct consequence of this construction we get the following
result.

\begin{Theorem} \label{teoremaReduccion}$\IA$ is a complete intersection $\Longleftrightarrow$ either
$\mA_{red} = \emptyset$ or $I_{\mA_{red}}$ is a complete
intersection.
\end{Theorem}

In \cite{BC}, the authors introduced the concept of free semigroup
to designate a family of subsemigroups of $\N$. Later, this concept
was generalized to subsemigroups of $\N^m$ in \cite{GR2}, where the
authors also proved several equivalent definitions for a semigroup
to be free. From the results here included one can derive that the
semigroup $\sum_{i = 1}^n \N a_i \subset \N^m$ is free if and only
if $\mA_{red} = \emptyset$.

Table \ref{algRed} shows an algorithm which receives as input the
set $\mA$ and computes $\mA_{red}$. The following example shows how
to compute $\mA_{red}$ following this algorithm.

\bigskip

\begin{table}[!htb]
\centering
\begin{tabular}{|p{8cm}|}
\hline
$$\begin{array}{l} \mbox{\bf Computation of } \mA_{red}\end{array}$$
$$\begin{array}{cl}
\ \mbox{Input:} & \mA = \{a_1,\ldots,a_n\} \subset \N^m \\
\ \mbox{Output:} & \mA_{red} \\
\end{array}
$$

\medskip
{
\begin{algorithmic}

\REPEAT

\STATE $\mB := \mA$

\FORALL {$a \in \mA$}

\IF {$a \notin \Q\, (\mA \setminus \{a\})$}

\STATE $\mA := \mA \setminus \{a\}$

\ELSE

\STATE $B := {\rm min}\{ b \in \Z^+\, \vert \, b\,a \in \Z\, (\mA
\setminus \{a\})\}$

\STATE $\mA := (\mA \setminus \{a\}) \cup \{B\, a\}$

 \IF {$B\, a \in \N (\mA \setminus \{B\, a\})$}

\STATE $\mA := \mA \setminus \{B\, a\}$

\ENDIF

\ENDIF

\ENDFOR

\UNTIL $(\mA = \emptyset)$ OR $(\mA = \mB)$

\RETURN ($\mA$)

\end{algorithmic}
}

\\

\hline
\end{tabular}

\medskip

\caption{Pseudo-code for computing $\mA_{red}$} \label{algRed}
\end{table}

\begin{Example}\label{ejemplored}
Set $\mA := \{a_1,a_2,a_3,a_4,a_5\} \subset \N^3$ with $a_1 =
(0,0,3),\, a_2 = (2,3,12),\, a_3 = (0,6,18),\, a_4 = (1,0,0)$ and
$a_5 = (1,5,17)$.

We firstly observe that $a_i \in \sum_{1 \leq j \leq 5 \atop j \neq
i} \Q\, a_j$ for all $i \in \{1,2,3,4,5\}$. We observe that $B_i =
1$ and that $B_i a_i = a_i \notin \sum_{j \in \{1,2,3,4,5\} \atop j
\neq i} \N a_j$ for all $i \in \{1,2,3,4\}$.

We compute $B_5$ and get that $B_5 = 3$, hence we denote $a_6 := 3
a_5$ and $\mA_1 := (\mA \setminus \{a_5\}) \cup \{a_6\} =
\{a_1,a_2,a_3,a_4,a_6\}$. Now we observe that $a_6 = a_1 + a_2 +
2a_3 + a_4 \in \N \{a_1,a_2,a_3,a_4\}$ and we write $\mA_2 := \mA_1
\setminus \{a_6\} = \{a_1,a_2,a_3,a_4\}$. We have that for all $i
\in \{1,2,3,4\}, a_i \in \sum_{j \in \{1,2,3,4\}\atop j \neq i} \Q\,
a_j$.

We compute $B_1^{\,\prime} := {\rm min}\{b \in \Z^+ \, \vert \, b
a_1 \in \sum_{j \in \{2,3,4\}} \Z\, a_j\}$ and we get that
$B_1^{\,\prime} = 2$. We denote $a_7 := 2 a_1$ and $\mA_3 := (\mA_2
\setminus \{a_1\}) \cup \{a_7\} = \{a_2,a_3,a_4,a_7\}$ and observe
that $a_7 \not\in \N \{a_2,a_3,a_4\}$.

We compute $B_2^{\,\prime} := {\rm min}\{b \in \Z^+ \, \vert \, b\,
a_2 \in \sum_{j \in \{3,4,7\}} \Z\, a_j\}$ and we get that
$B_2^{\,\prime} = 2$. We denote $a_8 := 2 a_2$ and $\mA_4 := (\mA_3
\setminus \{a_2\}) \cup \{a_8\} = \{a_3,a_4,a_7,a_8\}$. Finally we
observe that $a_8 = a_3 + 4a_4 + a_7 \in \N \{a_3,a_4,a_7\}$, and
write $\mA_5 := \mA_4 \setminus \{a_8\} = \{a_3,a_4,a_7\}$. Since
$a_3,a_4$ and $a_7$ are $\Q$-linearly independent, we get that
$\mA_{red} = \emptyset$.
\end{Example}

\medskip

\begin{Remark}\label{generadores} Whenever we
know by any means a set of $\mA$-homogeneous generators of
$I_{\mA_{red}}$ or $\mA_{red} = \emptyset$,  part (a) in {\rm Lemmas
\ref{redpertenece}, \ref{red2}} and {\rm Proposition \ref{red}} show
how to get a set of generators for $\IA$ without performing any
extra calculations. Moreover, this set of generators is minimal if
$\mA_{red} = \emptyset$ or the one of $I_{\mA_{red}}$ was. More
precisely, the set of generators of $\IA$ is obtained by following
these steps:

\begin{itemize}
\item[{\rm 1.}] If
$a_i = \sum_{j \in \{1,\ldots,n\} \atop j \neq i} \alpha_j a_j \in
\sum_{j \in \{1,\ldots,n\} \atop j \neq i} \N a_j$, then $\IA =
I_{\mA \setminus \{a_i\}} \cdot \kx + (x_i - \prod_{j \neq i}
x_j^{\alpha_j})$. Thus, if $\mathfrak B$ is a set of generators of
$I_{\mA \setminus \{a_i\}}$, then $\mathfrak B \cup \{x_i - \prod_{j
\neq i} x_j^{\alpha_j}\}$ is a set of generators of $\IA$. This is a
consequence of part {\rm (a)} in {\rm Lema \ref{redpertenece}}.

\item[{\rm 2.}] If $a_i
\notin \sum_{j \in \{1,\ldots,n\} \atop j \neq i} \Q\, a_j$, then
$\IA = I_{\mA \setminus \{a_i\}} \cdot \kx$. Thus, if $\mathfrak B$
is a set of generators of $I_{\mA \setminus \{a_i\}}$, then so is of
 $\IA$. This is a consequence of part {\rm (a)} in {\rm Lema \ref{red2}}.

\item[{\rm 3.}] If $a_i \in \sum_{j \in \{1,\ldots,n\} \atop j \neq i} \Q\, a_j$,
then denoting $\mA' := \{a_1,\ldots,a_{i-1},B_i a_i, a_{i+1},\ldots,
a_n\}$ we have that $\IA = \rho(I_{\mA'}) \cdot \kx$, where $\rho:
\kx \longrightarrow \kx$ is the $k$-homomorphism defined by $x_i
\mapsto x_i^{B_i}$ and $x_j \mapsto x_j$ if $j \neq i$. Thus, if
$\mathfrak B$ is a set of generators of $I_{\mA'}$, then $\mathfrak
\rho(\mathfrak B)$ is a set of generators of $\IA$. This is a
consequence of part {\rm (a)} in {\rm Proposition \ref{red}}.
\end{itemize}
\end{Remark}

Let us show how to get a minimal set of generators of the toric
ideal in {\rm Example \ref{ejemplored}} following the instructions
of {\rm Remark \ref{generadores}}.

\begin{Example}
\begin{itemize} \item[{\rm 1.}] We set $a_6 := 3 a_5$, $\mA_1 :=
\{a_1,a_2,a_3,a_4,a_6\}$ and denoting $$\rho_1:
k[x_1,x_2,x_3,x_4,x_6] \rightarrow k[x_1,x_2,x_3,x_4,x_5]$$ the
$k$-homomorphism defined by $\rho_1(x_i) = x_i$ for all $i \in
\{1,2,3,4\}$ and $\rho_1(x_6) = x_5^3$, then $\IA =
\rho_1(I_{\mA_1}) \cdot k[x_1,x_2,x_3,x_4,x_5]$.

\item[{\rm 2.}] We checked if $a_6 \in  \N a_1 + \N a_2 + \N a_3 + \N a_4$ and got that
 $a_6 = a_1 + 2a_2 + 2a_3 + a_4$. Hence, we set
$\mA_2 := \{a_1,a_2,a_3,a_4\}$, $g_1 := x_6 - x_1 x_2^2 x_3^2 x_4$
and have that $I_{\mA_1} = I_{\mA_2} \cdot k[x_1,x_2,x_3,x_4,x_6] +
(g_1).$

\item[{\rm 3.}] We defined $a_7 := 2 a_1$, $\mA_3 := \{a_2,a_3,a_4,a_7\}$ and denoting
$$\rho_2: k[x_2,x_3,x_4,x_7] \rightarrow
k[x_1,x_2,x_3,x_4]$$ the $k$-homomorphism defined by $\rho_2(x_i) =
x_i$ for all $i \in \{2,3,4\}$ and $\rho_2(x_7) = x_1^2$, then
$I_{\mA_2} = \rho_2(I_{\mA_3}) \cdot k[x_1,x_2,x_3,x_4]$.

\item[{\rm 4.}] We defined $a_8 := 2 a_2$, $\mA_4 := \{a_3,a_4,a_7,a_8\}$ and
denoting by  $$\rho_3: k[x_3,x_4,x_7,x_8] \rightarrow
k[x_2,x_3,x_4,x_7]$$ the $k$-homomorphism defined by $\rho_3(x_i) =
x_i$ for all $i \in \{3,4,7\}$ and $\rho_3(x_8) = x_2^2$, then
$I_{\mA_3} = \rho_3(I_{\mA_4}) \cdot k[x_2,x_3,x_4,x_7]$.

\item[{\rm 5.}] We checked if $a_8 \in \N a_3 + \N a_4 + \N a_7$ and got that
$a_8 = a_3 + 4 a_4 + a_7$. Hence, we set $\mA_5 := \{a_3,a_4,a_7\}$
and $g_2 := x_8 - x_3 x_4^4 x_7$ and have that $I_{\mA_4} =
I_{\mA_5}\cdot k[x_3,x_4,x_7,x_8] + (g_2)$.

\item[{\rm 6.}] Since $a_3, a_4, a_7$ are linearly independent, we have that $I_{\mA_5} = (0)$.
\end{itemize}

We finally obtained that $\mA_{red} = \emptyset$, therefore we
deduce that $\IA$ is a complete intersection and it is minimally
generated by the set of binomials $\{\rho_1(g_1),\,\rho_1 \circ
\rho_2 \circ \rho_3(g_2)\} = \{x_5^3 - x_1 x_2^2 x_3^2 x_4,\, x_2^2
- x_1^2 x_3 x_4^4\}$.
\end{Example}

\section{From $\IA$ to $I_{\mA_{(i,j)}}$}\label{sec3}
The objective of this section is to prove Theorem \ref{ci-iguales}.
This result provides, under certain hypotheses, necessary conditions
for $\IA$ to be a complete intersection. More precisely, for certain
$i,j: 1 \leq i < j \leq n$, we associate to $\IA$ a new toric ideal
$I_{\mA_{(i,j)}}$ in a ring of polynomials with exactly one variable
less and whose height is one unit less, such that $I_{\mA_{(i,j)}}$
is a complete intersection whenever $\IA$ is.

We denote by $\mH$ the set of elements of $\mA$ that belong to the
cone spanned by the rest of elements of $\mA$, i.e., $\mH :=
\left\{a_i \in \mA\, \vert \, a_i \in \sum_{j \in \{1,\ldots,n\}
\atop j \neq i} \R_{\geq 0}\, a_j \right\}.$ For all $a_i \in \mH$,
we set
\begin{center} $m_i := {\rm min}\left\{b \in \Z^+ \, \vert \,
b a_i \in \sum_{j \in \{1,\ldots,n\}  \atop j \neq i} \N
a_j\right\}$
\end{center}
and a binomial
\begin{center}$x_i^{m_i} - \prod_{j \in
\{1,\ldots,n\} \atop j \neq i} x_j^{\beta_{ij}} \in \IA$
\end{center}
is called a {\it critical binomial with respect to $x_i$}.

The concept of critical binomial was introduced by \cite{E} in the
context of toric ideals associated to affine monomial curves and
later studied by \cite{AV} in the same context. The definition
provided here is a natural extension of this concept to any toric
ideal. Critical binomials play an important role in the proofs of
the main results of this section and Section \ref{sec4}. It is well
known that critical binomials satisfy the following properties.

\begin{Lemma}\label{Lemacritico} \

\begin{itemize}
\item[{\rm (a)}] If $a_i \in \mH$ and $\mathfrak B$ is a set of generators of $\IA$
formed by binomials. Then there exists $g \in \mathfrak B$ such that
$g$ {\rm (}or $-g${\rm)} is a critical binomial with respect to
$x_i$.

\item[{\rm (b)}] If $f_1,\ldots,f_t$ are critical binomials with respect to
$x_{i_1},\ldots,x_{i_t}$ with $1 \leq i_1 < \cdots < i_t \leq n$ and
they all have different $\mA$-degrees, then the set
$\{f_1,\ldots,f_t\}$ can be extended to a minimal set of generators
of $\IA$ formed by binomials.
\end{itemize}
\end{Lemma}
%\begin{demo}Let $\mathfrak B = \{g_1,\ldots,g_{s}\}$ be a minimal set of generators of
%$\IA$ formed by binomials. We denote by $D_1,\ldots,D_{s} \in \N^m$
%the $\mA$-degrees of $g_1,\ldots,g_{s}$ respectively. If $f \in \IA$
%is a critical binomial with respect to $x_{i}$, then ${\rm
%deg}_{\mA}(f) = m_i a_i$ and $f$ can be written as $f =
%q_1g_1+\cdots+q_{s}g_{s}$\, where $q_{k}$ is either the null
%polynomial or an $\mA$-homogeneous polynomial of $\mA$-degree $m_i
%a_i - D_k \in \N^m$ for all $k \in \{1,\ldots,s\}$. If we apply to
%both sides of the equality the evaluation morphism defined by $x_k
%\mapsto 0$ for all $k \neq i$, then there exists $i_0 \in
%\{1,\ldots,s\}$ such that $q_{i_0} \neq 0$ and the image of
%$g_{i_0}$ under this morphism equals $x_i^{d}$ with $d a_i =
%D_{i_0}$. We assume that $i_0 = s$, then \,$g_s = x_i^{d} -
%\x^{\alpha}$\, and $D_{s} = d a_i \in \Z^+ a_i \cap
%\sum_{j\in\{1,\ldots,n\} \atop j \neq i} \N\,a_j$. By the definition
%of $m_i$ we get the equality $m_i = d$ and ${\rm deg}_{\mA}(q_{s}) =
%0$. As a consequence, $g_s$ (or $-g_s)$ is a critical binomial with
%respect to $x_i$, which proves (a). Moreover, $q_s \in k$, hence
%$\mathfrak B' := \{g_1,\ldots,g_{s-1},f\}$ is also a minimal set of
%binomials generating  $\IA$. If we iterate this process we get (b).
%\QED
%\end{demo}

\medskip Given $b = (b_1,\ldots,b_m)$, $c = (c_1,\ldots,c_m) \in
\N^m$, we denote by $\gcd\{b,c\} \in \N^m$ the vector whose $i$-th
coordinate is  $0$ if $b_i = 0$ or $c_i = 0$, or $\gcd\{b_i, c_i\}$
otherwise. Now we can formulate Theorem \ref{ci-iguales}.

\begin{Theorem}\label{ci-iguales}If there exist $i,j: 1 \leq i < j \leq n$ such that
$a_i,a_j \in \mH$ and $m_i a_i = m_j a_j$ and $\IA$ is a complete
intersection, then:

\begin{itemize}
\item[{\rm (a)}]\label{n-1}
$I_{\mA_{(i,j)}}$ is also a complete intersection, where
$\mA_{(i,j)} := \{a_1^{\,\prime},\ldots,a_{n-1}^{\,\prime}\}$ with
$a_k^{\,\prime}:= a_k$ for all $k \in \{1,\ldots,j-1\} \setminus
\{i\}$, $a_k^{\,\prime} := a_{k+1}$ for all $k\in
\{j\,,\ldots,n-1\}$ and $a_i^{\,\prime} := \gcd\{a_i, a_j\}.$
\item[{\rm (b)}]\label{mas}
For all $k \in \{1,\ldots,n-1\}$ such that $a_k^{\,\prime} \in {\rm
Cone}( \mA_{(i,j)} \setminus \{a_k^{\,\prime}\})$, if we denote
$m_k^{\,\prime} := \min\big\{b \in \Z^+\, \vert \, b\,
a_k^{\,\prime} \in \sum_{i\in\{1,\ldots,n-1\}\atop l \neq k}
\N\,a_{l}^{\,\prime} \big\}$, we have that $m_k^{\,\prime} = m_k$ if
$k \in \{1,\ldots,j-1\} \setminus \{i\}$, $m_k^{\,\prime} = m_{k+1}$
if $k \in \{j\,,\ldots,n-1\}$ and $m_k^{\,\prime} a_k^{\,\prime} \in
\N a_i + \N a_j$ if $k = i$.
\end{itemize}
\end{Theorem}

For proving Theorem \ref{ci-iguales} we will use Theorem 1.1 in
\cite{HSh}, which is a reformulation of \cite[Theorem
2.9]{Fisher-Shapiro}, and a technical lemma. Before presenting this
result we first include some definitions.

\begin{Definition}Let $B$ be an integral matrix. $B$
is {\it mixed} if every row of $B$ contains a positive and a
negative entry. $B$ is called {\it dominating} if it does not
contain a square mixed submatrix. For every $1 \leq t \leq {\rm
rk}(B)$, $\Delta_t(B)$ denotes the greatest common divisor of all
the nonzero $t \times t$ minors of $B$.
\end{Definition}

 For a binomial $f = \x^{\alpha} - \x^{\beta} \in \kx$, we
denote $\widehat f := \alpha - \beta \in \Z^n$.

\begin{Theorem}(\cite[Theorem 2.9]{Fisher-Shapiro}, \cite[Theorem 1.1]{HSh})
\label{comprueba} Let $\{f_1,\ldots,f_h\}$ be a set of binomials of
$\IA$ with $h = {\rm ht}(\IA)$ such that $f_i  = \x^{\alpha_i} -
\x^{\beta_i} \in \IA$ with
 $\gcd \{\x^{\alpha_i}, \x^{\beta_i}\} = 1$ for all $i \in \{1,\ldots,h\}$.
Let $A$ denote the $h \times n$ matrix whose $i$-th row is
$\widehat{f_i}$ for all $i \in \{1,\ldots,h\}$, then

\begin{center} $\IA = (f_1,\ldots,f_h) \,
 \Longleftrightarrow \, A$ is dominating and $\Delta_h(A) =  1$.
 \end{center}
\end{Theorem}

 In the proof of Theorem \ref{ci-iguales} we also use the following technical lemma.

\begin{Lemma}\label{dominantecontenido1} Let $A = (a_{(i,j)}) \in \mathcal M_{h \times n}(\Z)$
be an integral matrix of rank $h$ such that

\begin{itemize}
\item $a_{(h,j)} = 0$ for all $j \in \{1,\ldots,n-2\}$, \item
$a_{(h,n-1)} a_{(h,n)} < 0$ and $\gcd\{a_{(h,n-1)}, a_{(h,n)}\} =
1$.
\end{itemize}
Consider the matrix $A^{\,\prime} = (b_{(i,j)}) \in \mathcal M_{h -
1 \times n-1}(\Z)$ defined by
\begin{itemize}
\item $b_{(i,j)} := a_{(i,j)}$ for all $i \in \{1,\ldots,h-1\}$ and $j \in
\{1,\ldots,n-2\}$
\item $b_{(i,n-1)} := a_{(h,n)} a_{(i,n-1)} - a_{(h,n-1)} a_{(i,n)}$ for all $i \in
\{1,\ldots,h-1\}.$
\end{itemize}

Then,
\begin{itemize}  \item[{\rm (a)}] $A$ is dominating
$\Longleftrightarrow$ $A^{\,\prime}$ is dominating and $a_{(i,n-1)}
a_{(i,n)} \geq 0$ for all $i: 1 \leq i \leq h-1$.
\item[{\rm (b)}] $\Delta_h(A) = \Delta_{h-1}(A')$. \end{itemize}
\end{Lemma}

\begin{demo}For every $B \in \mathcal M_{u \times v}(\Z)$, we denote
by $B[i_1,\ldots,i_k][j_1,\ldots,j_k]$ the $k \times k$ submatrix
with rows $1 \leq i_1 < \cdots < i_k \leq u$, columns $1 \leq j_1 <
\cdots < j_k \leq v$ and $k \leq {\rm min}\{u, v\}$. Whenever, $v
\geq u$ we also denote by $B\{j_1,\ldots,j_u\}$ with $1 \leq j_1 <
\cdots < j_u \leq v$ the $u \times u$ minor of $B$ with rows
$1,\ldots,u$ and columns $j_1,\ldots,j_u$.

Suppose that $A$ is dominating, then $a_{(i,n-1)} a_{(i,n)} \geq 0$
for all $i: 1 \leq i \leq h-1$ because the submatrix $A[i,h][n-1,n]$
is not mixed. By contradiction assume that $A^{\,\prime}$ is not
dominating. Then there exist $k \leq h-1$, $1 \leq i_1 < \cdots <
i_k \leq h-1$ and $1 \leq j_1 < \cdots < j_k \leq n-1$ such that
$A^{\,\prime}[i_1,\ldots,i_k][j_1,\ldots,j_k]$ is mixed. Hence
$A[i_1,\ldots,i_k][j_1,\ldots,j_k]$ is mixed if $j_k < n-1$, or
$A[i_1,\ldots,i_k,h][j_1,\ldots,j_k,n]$ is mixed if $j_k = n-1$,
which is a contradiction.

Suppose now that $A$ is not dominating and that $a_{(i,n-1)}
a_{(i,n)} \geq 0$ for all $i: 1 \leq i \leq h-1$, then there exist
$k \leq h$, $1 \leq i_1 < \cdots < i_k \leq h$ and $1 \leq j_1 <
\cdots < j_k \leq n$ such that $A[i_1,\ldots,i_k][j_1,\ldots,j_k]$
is mixed, then it is easy to check that
\begin{itemize}
\item $A^{\,\prime}[i_1,\ldots,i_k][j_1,\ldots,j_k]$ is mixed if $i_k < h$ and $j_k < n-1$;
\item $A^{\,\prime}[i_1,\ldots,i_k][j_1,\ldots,j_{k-1},n-1]$ is
mixed if $i_k < h$, $j_{k-1} < n-1$ and $j_k \geq n-1$; or
\item $A^{\,\prime}[i_1,\ldots,i_{k-1}][j_1,\ldots,j_{k-1}]$ is mixed otherwise.
\end{itemize}
Hence $A^{\,\prime}$ is not dominating and we have (a).

To prove (b) it suffices to observe that  $A\{j_1,\ldots,j_{h}\}$
equals:
\begin{itemize}
\item $0$, if $j_{h} < n-1$;
\item $\pm a_{(h,n-1)} A^{\,\prime}\{j_1,\ldots,j_{h-1}\}$, if $ j_{h} = n -
1$;
\item $\pm a_{(h,n)} A^{\,\prime}\{j_1,\ldots,j_{h-1}\}$, if $j_{h-1}
< n-1$ and $ j_{h} = n$; or
\item $\pm A^{\,\prime}\{j_1,\ldots,j_{h-1}\}$, if $j_{h-1} = n-1$ and $j_{h} =
n$,
\end{itemize}
and that $\gcd\{a_{(h,n-1)},a_{(h,n)}\} = 1$. Hence $\Delta_{h}(A) =
\Delta_{h-1}(A^{\,\prime})$. \QED
\end{demo}

\bigskip

{\noindent \bf Proof of Theorem \ref{ci-iguales}.} We can assume
without loss of generality that $i = n-1$ and $j = n$. We denote
$\mA' := \mA_{(n-1,n)} = \{a_1,\ldots,a_{n-2},a_{n-1}'\}$, where
$a_{n-1}' = \gcd\{a_{n-1},a_n\}$, and let us prove that $I_{\mA'}$
is a complete intersection. From the definition of $m_{n-1}$ and
$m_n$ it follows that $\gcd\{m_{n-1},m_n\} = 1$ and it is easy to
check that $m_n a_{n-1}' = a_{n-1}$ and $m_{n-1} a_{n-1}' = a_{n}$.
Since $\Q\, \mA = \Q\, \mA^{\,\prime}$, then ${\rm
ht}(I_{\mA^{\,\prime}}) = {\rm ht}(\IA) - 1$.

Consider the critical binomial $f := x_n^{m_n} - x_{n-1}^{m_{n-1}}$
with respect to $x_n$. By Lemma \ref{Lemacritico} (b), there exists
$\{f_1,\ldots,f_h\}$ a minimal set of generators of $\IA$ formed by
binomials such that $f_h = f$. On the other hand, $h = {\rm
ht}(\IA)$ because $\IA$ is a complete intersection.

Consider now  $\psi: \kx \longrightarrow k[x_1,\ldots,x_{n-1}]$, the
 $k$-homomorphism defined by $\psi(x_k) = x_k$ for all $k \in
\{1,\ldots,n-2\}$, $\psi(x_{n-1}) = x_{n-1}^{m_n}$ and $\psi(x_n) =
x_{n-1}^{m_{n-1}}$. Let us prove that
\begin{equation}\label{igualdad} I_{\mA^{\,\prime}} =
(\psi(f_1),\ldots,\psi(f_{h-1})),
\end{equation}
which implies that $I_{\mA^{\,\prime}}$ is a complete intersection
and the first part of the theorem holds. Indeed, for all $\alpha \in
\N^n$ we have that ${\rm deg}_{\mA}(\x^{\alpha}) = {\rm
deg}_{\mA^{\,\prime}}(\psi(\x^{\alpha}))$, which implies that
$\{\psi(f_1),\ldots,\psi(f_{h-1})\} \subset I_{\mA'}$. If we
consider $A$ the $h \times n$ integral matrix whose $k$-th row is
$\widehat{f_k} \in \Z^n$ for all $k \in \{1,\ldots,h\}$, from
Theorem \ref{comprueba} it follows that $A$ is dominating and
$\Delta_{h}(A) = 1$ because  $\IA = (f_1,\ldots,f_h)$. Then, if we
take $A'$  the $(h-1) \times (n-1)$ matrix whose $k$-th row is
$\widehat{\psi(f_k)} \in \Z^{n-1}$ for all $k \in \{1,\ldots,h-1\}$,
then $A^{\,\prime}$ is dominating and $\Delta_{h-1}(A^{\,\prime}) =
\Delta_h(A) = 1$ by Lemma \ref{dominantecontenido1}. Then, by virtue
of Theorem \ref{comprueba} we conclude that $I_{\mA^{\,\prime}}=
(\psi(f_1),\ldots,\psi(f_{h-1}))$.

In accordance with our previous assumptions, let us now prove the
second part of the theorem. Let us prove first that $m_k =
m_k^{\,\prime}$ for all $a_k \in {\rm Cone}(\mA' \setminus \{a_k\})$
with $k \in \{1,\ldots ,n-2\}$\,, where $m_k^{\,\prime}$ equals
$\min\left\{ b \in \Z^+ \, \vert \, b a_{k} \in \sum_{l
\in\{1,\ldots,n-2\} \atop l \neq k} \N\, a_l + \N\,
a_{n-1}^{\,\prime} \right\}$. Take $k \in \{1,\ldots,n-2\}$ such
that $a_k \in {\rm Cone}(\mA^{\,\prime} \setminus \{a_k\})$, the
inequality $m_k \geq m_k^{\,\prime}$\, is obvious because
$a_{n-1},\,a_n \in \N\, a_{n-1}^{\,\prime}$. By Lemma
\ref{Lemacritico} (a), there exists $l \in \{1,\ldots,h-1\}$ such
that $\psi(f_{l})$ is a critical binomial with respect to $x_k$.
Then \,$f_{l}=x_k^{m_k^{\,\prime}}-\x^{\alpha}$\,, where
$\x^{\alpha}$ is a monomial of $\mA$-degree $m_k^{\,\prime} a_k$ not
involving the variable $x_k$, and hence $m_k^{\,\prime} a_k \in
\sum_{t \in\{1,\ldots,n\} \atop t \neq k} \N\,a_t$. This implies
that $m_k^{\,\prime} \geq m_k$ by the definition of $m_k$, thus the
equality $m_k^{\,\prime} = m_k$ follows. Finally suppose that
 $a_{n-1}^{\,\prime} \in {\rm Cone}(\mA'
\setminus \{a_{n-1}^{\,\prime}\})$ and let us see that $m_{n-1}'
a_{n-1}' \in \N a_{n-1} + \N a_n$. Indeed, again by
 Lemma \ref{Lemacritico} (a), there exists $l \in \{1,\ldots,h-1\}$ such that
 $\psi(f_{l})$ is a critical binomial with respect to  $x_{n-1}$.
Then \,$f_{l}=x_{n-1}^{\beta_{n-1}} x_n^{\beta_n} - \x^{\alpha}$\,,
where $x_{n-1}^{\beta_{n-1}} x_n^{\beta_n}$ is a monomial of
$\mA$-degree $m_{n-1}^{\,\prime} a_{n-1}^{\,\prime}$, and hence
$m_{n-1}^{\,\prime} a_{n-1}^{\,\prime} \in \N\,a_{n-1} + \N\, a_n$.
\QED

\medskip

The following example shows that the necessary conditions for $\IA$
to be a complete intersection of Theorem \ref{ci-iguales} are not
sufficient in general.

\begin{Example}\label{NCI}\cite[Example 2.2]{BGS} The toric ideal $\IA$ with
$\mA := \{45,70,75,98,147\} \subset \N$ is not a complete
intersection because ${\rm ht}(\IA) = 4$ and $\mu(\IA) = 7$.
Nevertheless, we have that

\begin{center}
$m_1 := {\rm min}\{b \in \Z^+ \, \vert\, b a_1 \in \N a_2 + \N a_3 +
\N a_4 + \N a_5\} = 5,$ and

$m_3 := {\rm min}\{b \in \Z^+ \, \vert\, b a_3 \in \N a_1 + \N a_2 +
\N a_4 + \N a_5\} = 3,$
\end{center}
and we have that $m_1 a_1 = m_3 a_3 = 225$. Moreover, denoting
$a_1^{\,\prime} := \gcd\{a_1,a_3\}= 15$, $a_2^{\,\prime} := 70$,
$a_3^{\,\prime} := 98$ and $a_4^{\,\prime} := 147$ we have that
$\mA_{(1,3)} = \{a_1^{\,\prime}, a_2^{\,\prime}, a_3^{\,\prime},
a_4^{\,\prime}\}$ and $I_{\mA_{(1,3)}}$ is a complete intersection.
Finally, denoting

\begin{center} $m_i^{\,\prime} := {\rm min}\{b \in \Z^+ \, \vert \,
b a_i^{\,\prime} \in \sum_{j \in \{1,2,3,4\} \atop i \neq j} \N \,
a_j^{\,\prime} \}$ for all $i \in \{1,2,3,4\},$
\end{center} it is easy to check that  $m_1^{\,\prime} = 14$, $m_1^{\,\prime}
a_1^{\,\prime} \in \N a_1 + \N a_3$, $m_2^{\,\prime} = m_2 = 3$\,,
$m_3^{\,\prime} = m_4 = 3$ and $m_4^{\,\prime} = m_5 = 2\,.$
\end{Example}

\medskip We end this section with a technical lemma which is a
generalization of part (b) in Theorem \ref{ci-iguales}. Indeed,
Theorem \ref{ci-iguales}.(b) is obtained from Lemma
\ref{lematecnico} if we consider the unitary sets
 $V' = \{a_k'\}$ for all $k \in \{1,\ldots,n-1\}$ such that $a_k' \in {\rm Cone}(\mA_{(i,j)}
\setminus \{a_k'\})$. We will use this result in the next section to
prove correctness of Algorithm CI-simplicial.

\begin{Lemma}\label{lematecnico}Suppose that there exist $a_{i},\,a_{j} \in \mH$ with $1 \leq i < j \leq n$
such that $m_{i} a_{i} = m_j a_j$ and consider

\begin{itemize}
\item $V' \subsetneq \mA_{(i,j)}$ such that ${\rm
dim}(\Q V') = 1$ and  $a := \gcd(V') \in {\rm Cone}(\mA_{(i,j)}
\setminus V')$,
\item $M := {\rm min}\{b \in \Z^+ \, \vert \, b a \in \sum_{a_k' \in
\mA_{(i,j)} \setminus V'} \N a_k'\}$, and \item $V := V'$ if $a_i' =
\gcd\{a_i,a_j\} \notin V'$, or $V := (V' \setminus \{a_i'\}) \cup
\{a_{i},a_j\}$ otherwise.
\end{itemize}

 If $\IA$ is a complete intersection and $M a \in \N V'$,
then $M a \in \N V \cap \N (\mA \setminus V)$.
\end{Lemma}
\begin{demo}Suppose without loss of generality that $i = n-1$, $j =
n$ and set $\mA^{\,\prime} = \mA_{(n-1,n)} =
\{a_1',\ldots,a_{n-1}'\}$, where $a_i' = a_i$ for all $i \in
\{1,\ldots,n-2\}$ and $a_{n-1}' = \gcd\{a_{n-1},a_n\}$. By Theorem
\ref{ci-iguales} (a) we have that $I_{\mA^{\,\prime}} \subset
k[x_1,\ldots,x_{n-1}]$ is a complete intersection. Moreover, by
(\ref{igualdad}) in the proof of Theorem \ref{ci-iguales}, if we
denote by  $\psi$ the $k$-homomorphism $\psi: \kx \longrightarrow
k[x_1,\ldots,x_{n-1}]$ defined by $\psi(x_i) = x_i$ for $i \in
\{1,\ldots,n-2\}$, $\psi(x_{n-1}) = x_{n-1}^{m_n}$, $\psi(x_n) =
x_{n-1}^{m_{n-1}}$; then there exist $h = {\rm ht}(\IA)$ binomials
$g_1,\ldots,g_{h} \in \IA$ such that $I_{\mA^{\,\prime}} =
(\psi(g_1),\ldots,\psi(g_{h-1}))\subset k[x_1,\ldots,x_{n-1}].$

Since $M a \in \N V^{\,\prime} \cap \N (\mA' \setminus
V^{\,\prime})$, it follows that $M a = \sum_{a_u' \in V^{\,\prime}}
\gamma_u a_u' = \sum_{a_v' \in \mA^{\,\prime} \setminus
V^{\,\prime}} \gamma_v a_v'$. Then, denoting $\alpha := \sum_{a_u'
\in V^{\,\prime}} \gamma_u \bar{e}_u$ and $\beta := \sum_{a_v' \in
\mA^{\,\prime} \setminus V^{\,\prime}} \gamma_v \bar{e}_v$, where
$\{\bar{e}_1,\ldots,\bar{e}_{n-1}\}$ is the canonical basis of
$\Z^{n-1}$, we have the binomial $f := \x^{\alpha} - \x^{\beta} \in
I_{\mA^{\,\prime}}$ whose $\mA'$-degree is $M a \in \N^m$.
Consequently, $f = q_1 \psi(g_1) + \cdots + q_{h-1} \psi(g_{h-1})$,
where for all $t \in \{1,\ldots,h-1\}$ $q_t$ is either the null
polynomial or an $\mA^{\,\prime}$-homogeneous polynomial of
$\mA^{\,\prime}$-degree $M a - {\rm
deg}_{\mA^{\,\prime}}(\psi(g_{t})) \in \N^m$.

Consider the $k$-homomorphism $\Phi: k[x_1,\ldots,x_{n-1}]
\longrightarrow k[x_1,\ldots,x_{n-1}]$ defined by $x_k \mapsto 0$
for all $k \in \{1,\ldots,n-1\}$ such that $a_k' \notin V'$. Then,
$\x^{\alpha} = \Phi(q_1) \Phi(\psi(g_1)) + \cdots + \Phi(q_{h-1})
\Phi(\psi(g_{h-1}))$ and there exist
$\lambda_1,\ldots,\lambda_{n-1},\delta_1,\ldots \delta_{n-1} \in \N$
y $i_0 \in \{1,\ldots,h-1\}$ verifying the following:
\begin{itemize}
\item
 $\Phi(q_{i_0})\neq 0$,

\item
$\psi(g_{i_0}) = x_1^{\lambda_1} \cdots x_{n-1}^{\lambda_{n-1}}
 - x_1^{\delta_1} \cdots x_{n-1}^{\delta_{n-1}}$,

\item
$\lambda_r=0$ if $a_r' \notin V'$ for all  $r \in \{1,\ldots,n-1\}$,
and

\item
$\delta_s \neq 0$ for some $s \in \{1,\ldots,n-1\}$ such that $a_s'
\notin V'$.
\end{itemize}

Then, $M a$ is componentwise $\geq$ than $\lambda_1 a_1' + \cdots +
\lambda_{n-1} a_{n-1}' = \delta_1 a_1' + \cdots +
 \delta_{n-1} a_{n-1}'$, and $M a - \sum_{a_u' \in
V^{\,\prime}} \delta_u a_u'$ is componentwise $\geq$ than
 $\lambda_1 a_1' + \cdots + \lambda_{n-1} a_{n-1}' - \sum_{a_u' \in
V^{\,\prime}} \delta_u  a_u' =  \sum_{a_u' \in \mA^{\,\prime}
\setminus V^{\,\prime}} \delta_u a_u'$, which belongs to
$$\Z^+ a \, \bigcap\, \sum_{a_u' \in \mA^{\,\prime}
\setminus V^{\,\prime}} \N\,a_u'.$$ From the definition of $M$ we
get that $\delta_u = 0$ for all $a_u' \in V'$ and $M a = {\rm
deg}_{\mA^{\,\prime}}(\psi(g_{i_0})) = {\rm deg}_{\mA}(g_{i_0})$.

If $a_{n-1}' \in V'$, then $g_{i_0} = x_{1}^{\lambda_1} \cdots
x_{n-2}^{\lambda_{n-2}} x_{n-1}^{\overline{\lambda}_{n-1}}
x_n^{\overline{\lambda}_n} - x_1^{\delta_1} \cdots
x_{n-2}^{\delta_{n-2}}$, where
$\overline{\lambda}_{n-1},\overline{\lambda}_{n}$ are nonnegative
integers such that $\overline{\lambda}_{n-1} a_{n-1}
\,+\,\overline{\lambda}_{n} a_n=\lambda_{n-1} a_{n-1}'$. As a
consequence, $M a
 = \lambda_1 a_1 + \cdots + \lambda_{n-2} a_{n-2} + \overline{\lambda}_{n-1}
a_{n-1} + \overline{\lambda}_{n} a_n \in \N V$.

If $a_{n-1}' \notin V'$, then $g_{i_0} = x_{1}^{\lambda_1} \cdots
x_{n-2}^{\lambda_{n-2}} - x_1^{\delta_1} \cdots
x_{n-2}^{\delta_{n-2}}  x_{n-1}^{\overline{\delta}_{n-1}}
x_n^{\overline{\delta}_{n}}$, where $\overline{\delta}_{n-1},
\overline{\delta}_{n}$ are nonnegative integers such that
$\overline{\delta}_{n-1} a_{n-1} \,+\,\overline{\delta}_{n} a_n =
\delta_{n-1} a_{n-1}'$. As a consequence, $M a
 = \delta_1 a_1 + \cdots + \delta_{n-2} a_{n-2} + \overline{\delta}_{n-1}
a_{n-1} + \overline{\delta}_{n} a_n \in \N (\mA \setminus V)$. \QED
\end{demo}

\section{Complete intersection simplicial toric ideals} \label{sec4}

This section concerns simplicial toric ideals, it is worth
mentioning that given a set $\mA = \{a_1,\ldots,a_n\} \subset \N^m$,
one can check if $\IA$ is a simplicial toric ideal by obtaining a
minimal set of generators $\{a_{i_1},\ldots,a_{i_s}\}$ of ${\rm
Cone}(\mA)$ and checking if $s = {\rm dim}(\Q \mA)$. One could
obtain a minimal set of generators ${\rm Cone}(\mA)$ if one knows
how to determine whether $a_i \in {\rm Cone}(\mA \setminus
\{a_i\})$. Indeed, this problem is equivalent to check if the system
of equations $ \sum_{j \in \{1,\ldots,n\} \atop j \neq i} x_j a_j =
a_i$, with $x_j \geq 0$ for all $j \in \{1,\ldots,n\} \setminus
\{i\}$ is feasible; this can be done by means of the simplex method,
see, e.g., \cite[Section 5.2]{Dantzig}.

In this section we study the property of being a complete
intersection in simplicial toric ideals. More precisely, we focus on
the design of Algorithm CI-simplicial, an algorithm for checking if
a simplicial toric ideal is a complete intersection. This algorithm
arises as a consequence of the convenient application of Theorems
\ref{teoremaReduccion} and \ref{ci-iguales} and Lemma
\ref{lematecnico} to the simplicial context together with
Proposition \ref{simplicial}, which is a specific result for
simplicial toric ideals.

\begin{Proposition}\label{simplicial}Let $\IA$ be a simplicial toric
ideal. If $\IA$ is a complete intersection, then one of the
following holds:
\begin{itemize}
\item[{\rm (a)}] there exist $i,j: 1 \leq i < j \leq n$ such that $a_i, a_j \in \mH$ and $m_i a_i = m_j a_j$, or
\item[{\rm (b)}] $\mA_{red} = \emptyset$.
\end{itemize}
\end{Proposition}
\begin{demo}Denote $r := {\rm dim}(\Q \mA)$ and suppose that $m_i a_i \neq m_j a_j$
for every $a_i, a_j \in \mH$, $i \neq j$. We may assume without loss
of generality that $\{a_1,\ldots,a_r\}$ is a minimal set of
generators of the cone spanned by $\mA$, then
$\{a_{r+1},\ldots,a_n\} \subset \mH$ and whenever $i,j: r+1 \leq i <
j \leq n$ we have that $m_i a_i \neq m_j a_j$. We aim prove that
$\mA_{red} = \emptyset$. Let $f_{r+1},\ldots,f_n$ be critical
binomials with respect to $x_{r+1},\ldots,x_n$ respectively. Since
${\rm deg}_{\mA}(f_i) = m_i a_i \neq m_j a_j = {\rm deg}_{\mA}(f_j)$
for every $r+1 \leq i < j \leq n$, by Lemma \ref{Lemacritico}.(b) we
know that there exists a minimal set of generators of $\IA$
containing $\{f_{r+1},\ldots,f_n\}$. Moreover, ${\rm ht}(\IA) = n-r$
and $\IA$ is a complete intersection, thus $\IA =
(f_{r+1},\ldots,f_n)$.

We write $f_i = x_i^{m_i} - \x^{\alpha_i}$ for every $i: r+1 \leq i
\leq n$ and we claim that there exists $j \in \{r+1,\ldots,n\}$ such
that $x_j \nmid \x^{\alpha_k}$ for every $k \in \{r+1,\ldots,n\}$.
Assume this claim is false and consider the simple directed graph
with vertex set $\{r+1,\ldots,n\}$ and arc set $\{(j,k)\, \vert\, r
+ 1\leq j,\,k \leq n$ and $x_j \mid \x^{\alpha_k}\}$, it is clear
that the out-degree of every vertex is greater or equal to one,
which implies that there is a cycle in the graph. Suppose, without
loss of generality, that the cycle is $(r+1,r+2,\ldots,r+k,r+1)$
with $k \leq n-r$, this means that $(f_{r+1},\ldots,f_{r+k}) \subset
(x_{r+1},\ldots,x_{r+k})$, so $\IA \subsetneq H :=
(x_{r+1},\ldots,x_{r+k},f_{r+k+1},\ldots,f_{n})$ but this is not
possible because $\IA$ is prime and $n-r = {\rm ht}(\IA) < {\rm
ht}(H) \leq n-r$.

Thus there exists $i \in \{r+1,\ldots,n\}$ such that $x_i \nmid
\x^{\alpha_j}$ for every $j \in \{r+1,\ldots,n\}$. Suppose that $i =
n$ and let us prove that $B_n a_n \in \sum_{j = 1}^{n-1} \N a_j$. By
\cite[Proposition 2.3]{EV},
$\{\widehat{f_{r+1}},\ldots,\widehat{f_{n}}\}$ is a $\Z$-basis for
the kernel of the homomorphism $\pi: \Z^n \longrightarrow \Z^m$
induced by $\pi(e_j) = a_j$. By definition, $B_n a_n = \sum_{j =
1}^{n-1} \beta_j a_j$ for some $\beta_1,\ldots,\beta_{n-1} \in \Z$,
so take $\delta := B_n e_n - \sum_{j = 1}^{n-1} \beta_j e_j$, then
$\delta \in {\rm ker}(\pi)$. Consequently, if we express $\delta$ as
a combination of $\widehat{f_{r+1}},\ldots,\widehat{f_{n}}$ we
derive that $m_n
 \mid B_n$, so $B_n a_n
\in \sum_{j = 1}^{n-1} \N a_j$. Now, by Proposition \ref{red} and
Lemma \ref{redpertenece}, it follows that $I_{\mA \setminus
\{a_n\}}$ is a complete intersection minimally generated by
$\{f_{r+1},\ldots,f_{n-1}\}$. Moreover $f_i \in I_{\mA \setminus
\{a_n\} }$ is a critical binomial with respect to $x_i$ for all $i
\in \{r+1,\ldots,n-1\}$. If we iterate the same argument we get that
we can reorder $a_{r+1},\ldots,a_{n-1}$ in such a way that $B_i a_i
\in \sum_{j = 1}^{i-1} \N a_j$ for all $i \in \{r+1,\ldots,n-1\}$.
Then by Proposition \ref{red} and Lemma \ref{redpertenece} we obtain
that $\mA_{red} = \mB_{red}$, where $\mB = \{a_1,\ldots,a_r\}$.
Since $a_1,\ldots,a_r$ are $\Q$-linearly independent, by Lemma
\ref{red2} we deduce that $\mB_{red} = \emptyset$ and the proof is
complete. \QED
\end{demo}

\medskip
The condition of being a simplicial toric ideal in this proposition
is essential to obtain the result because there exist complete
intersection toric ideals such that $\mA_{red} \neq \emptyset$ and
$m_i a_i \neq m_j a_j$ for all $a_i, a_j \in \mH$. Let us see an
example.

\begin{Example}Let $\IA$ be the toric ideal associated to $\mA =
\{a_1,a_2,a_3,a_4,a_5\} \subset \N^3$, where $a_1 = (0,2,1),\,a_2 =
(4,2,1),\, a_3 = (2,2,1),\, a_4 = (1,3,1)$ and $a_5 = (1,1,1)$. This
toric ideal is not simplicial because ${\rm dim}(\Q \mA) = 3$ and
${\rm Cone}(\mA)$ has $4$ extremal rays; indeed
$\{a_1,a_2,a_4,a_5\}$ is a minimal set of generators of ${\rm
Cone}(\mA)$.

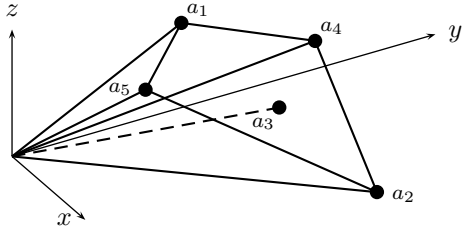
\psset{Alpha=120,unit=1.3}
\begin{pspicture}(-3,0)(2,2)

\pstThreeDCoor[linecolor = black,xMin=0,xMax=1.5,yMin=0,yMax=5,
zMin=0,zMax=1.5]

\psset{arrowscale=1.5,arrowinset=0,dotstyle=*,dotscale=1.5}
\pstThreeDLine(0,0,0)(0,2,1) \pstThreeDLine(0,0,0)(4,2,1)
\pstThreeDLine(0,0,0)(1,3,1) \pstThreeDLine(0,0,0)(1,1,1)
\pstThreeDLine[linestyle = dashed](0,0,0)(2,2,1)

\pstThreeDLine(0,2,1)(1,3,1) \pstThreeDLine(0,2,1)(1,1,1)
\pstThreeDLine(4,2,1)(1,3,1) \pstThreeDLine(4,2,1)(1,1,1)

\pstThreeDDot(0,2,1) \pstThreeDDot(4,2,1) \pstThreeDDot(2,2,1)
\pstThreeDDot(1,3,1) \pstThreeDDot(1,1,1)

\pstThreeDPut(0,2.2,1.1){\footnotesize $a_1$}
\pstThreeDPut(4.2,2.2,1){\footnotesize $a_2$}
\pstThreeDPut(2.2,1.7,1){\footnotesize $a_3$}
\pstThreeDPut(1,3.2,1.1){\footnotesize $a_4$}
\pstThreeDPut(1,.7,1.1){\footnotesize $a_5$}
\end{pspicture}

\vspace{2cm} $\IA$ is a height $2$ ideal generated by $g_1 := x_3^2
- x_1 x_2$ and $g_2 := x_1 x_3 - x_4 x_5$, thus it is a complete
intersection. Nevertheless, $\mH = \{a_3\}$. Moreover, the relations
$2\,a_3 = a_1 + a_2$ and $a_1 + a_3 = a_4 + a_5$ show that $B_i = 1$
for all $i \in \{1,2,3,4,5\}$ and $\mA$ is a minimal set of
generators of the semigroup $\sum_{i = 1}^5 \N a_i$, thus $\mA_{red}
= \mA$.\end{Example}

In Table \ref{algoritmo} we propose Algorithm CI-simplicial, which
works in the following way. It receives as input a set $\mA \subset
\N^m$ such that $\IA \subset k[x_1,\ldots,x_n]$ is a simplicial
toric ideal. If there exist $a_i, a_j \in \mathcal H$ such that $m_i
a_i = m_j a_j$, we consider the set $\mA_{(i,j)} = (\mA \setminus
\{a_i,a_j\}) \cup \{\gcd\{a_i,a_j\}\} \subset \N^m$.

 The following properties hold for $I_{\mA_{(i,j)}}$:

\begin{itemize} \item[(1)] $I_{\mA_{(i,j)}} \subset k[x_1,\ldots,x_{n-1}]$
and ${\rm ht}(I_{\mA_{(i,j)}}) = {\rm ht}(\IA) - 1$,
\item[(2)] $I_{\mA_{(i,j)}}$ is a complete intersection whenever $\IA$ is
(see Theorem \ref{ci-iguales} (a)), and
\item[(3)] $I_{\mA_{(i,j)}}$ is a simplicial toric ideal (because ${\rm Cone}(\mA)
= {\rm Cone}(\mA_{(i,j)})$).
\end{itemize}

 Proceeding as in Theorem \ref{ci-iguales} (b), if $a_i' \in
{\rm Cone}(\mA_{(i,j)} \setminus \{a_i'\})$ and $m_i' a_i' \notin \N
a_i + \N a_j$, where $a_i' = \gcd\{a_i,a_j\}$, then $\IA$ is not a
complete intersection. Otherwise, we iterate this procedure as many
times as possible until we get a set $\mB = \{b_1,\ldots,b_{n'}\}
\subset \N^m$ such that $\IB$ is a simplicial toric ideal satisfying
that if $\bar{\mathcal H} := \{b_i \in \mB \, \vert \, b_i \in {\rm
Cone}(\mB \setminus \{b_i\}\}$ and $\bar{m}_i := {\rm min}\{c \in
\Z^+ \, \vert \, c b_i \in \sum_{j \in \{1,\ldots,n'\} \atop j \neq
i} \N b_j\}$ for all $b_i \in \bar{\mathcal H}$, then $\bar{m}_i b_i
\neq \bar{m}_j b_j$ for all $b_i, b_j \in \bar{\mathcal H}$.

 Then we compute $\mB_{red}$. If $\mB_{red} \not=
\emptyset$, we conclude that $\IB$ is not a complete intersection by
Proposition \ref{simplicial}, hence $\IA$ is not a complete
intersection by Theorem \ref{ci-iguales} (a) and we are done. In
case $\mB_{red} = \emptyset$, then $\IB$ is a complete intersection
by Theorem \ref{teoremaReduccion}, but we can not assert whether
$\IA$ is a complete intersection. In order to decide if $\IA$ is a
complete intersection we have to check several additional conditions
consisting in determining whether certain elements of $\N^m$ belong
to some semigroups. If any of these elements does not belong to its
corresponding semigroup, then by Lemma \ref{lematecnico} we get that
$\IA$ is not a complete intersection. Otherwise, we can ensure that
$\IA$ is a complete intersection.

\begin{table}
\centering
\begin{tabular}{|p{12cm}|}
\hline
\begin{center}
{\bf Algorithm~CI-simplicial}
\end{center}

$$\begin{array}{cl}
\ \mbox{Input:} & \mA = \{a_1,\ldots,a_n\} \subset \N^m \mbox{ such
that } \IA \mbox{ is a simplicial toric ideal}
\\
\ \mbox{Output:} & \mbox{{\sc True} or {\sc False}}  \\
\end{array}
$$

\medskip

\begin{algorithmic}

\STATE $G := \mA$

\STATE $V_i := \{a_i\}$, $\forall i \in \{1,\ldots,n\}$

\STATE $m_{i} := \min\{b \in \Z^+\, \vert\, b  a_{i} \in \N\,(\mA
\setminus \{a_{i}\})\}$ for all $a_i \in {\rm Cone}(\mA \setminus
\{a_i\})$

\STATE $k := 0$

\WHILE {$\exists\, a_i, a_j \in G$ such that $m_i a_i = m_j a_j$}

\IF {$m_{i} a_{i} \not\in \N\, V_{i} \cap \N \, V_j$}

\RETURN {\sc False}

\ENDIF

\STATE $k := k + 1$; $a_{n+k} :=\gcd\{a_i,\, a_j\}$; $G := (G
\setminus \{ a_i,\, a_j \}) \cup \{a_{n+k}\}$; $V_{n+k} := V_i \cup
V_j$

\IF {$a_{n+k} \in {\rm Cone}(G \setminus \{a_{n+k}\})$}

\STATE $m_{n+k} :=\min\{b \in \Z^+\, \vert \, b\, a_{n+k} \in
 \N\, G \}$

\ENDIF

\ENDWHILE

$\mB := G$

\REPEAT

\STATE $G := \mB$

\FORALL {$a_i \in G$}

\IF {$a_i \notin \Q\, (\mB \setminus \{a_i\})$}

\STATE $\mB := \mB \setminus \{a_i\}$

\ELSE

\STATE $B_i := {\rm min}\{ b \in \Z^+\, \vert \, b a_i \in \Z\,(\mB
\setminus \{a_i\}) \}$

\IF {$B_i a_i \in \N (\mB \setminus \{a_i\})$}

\IF {$B_i a_i \notin \N V_i \cap \sum_{a_j \in \mB \setminus
\{a_i\}} \N\, V_j $}

\RETURN {\sc False}

\ENDIF

\STATE  $\mB := \mB \setminus \{a_i\}$

\ENDIF

\ENDIF

\ENDFOR

\UNTIL $(\mB = \emptyset)$ OR $(\mB = G)$

\IF {$\mB \not= \emptyset$}

\RETURN {\sc False}

\ENDIF

\RETURN {\sc True}

\end{algorithmic}

\\
\hline
\end{tabular}

\medskip

 \caption{Pseudo-code for checking whether a simplicial
toric ideal is a complete intersection} \label{algoritmo}
\end{table}

\medskip Our next goal is to prove correctness of this algorithm.

\begin{Theorem}\label{algoritmosimplicial}
Let $\IA$ be a simplicial toric ideal,  {\rm
Algorithm~CI-simplicial} determines whether $\IA$ is a complete
intersection.
\end{Theorem}
\begin{demo} The algorithm always terminates. Moreover, if
$m_i a_i \neq m_j a_j$ whenever $a_i, a_j \in \mH$, then the
algorithm does not enter the {\bf while} loop. Thus one can observe
that {\rm CI-simplicial}$(\mA) =$ {\sc True} if and only if
$\mA_{red} = \emptyset$. By Proposition \ref{simplicial} and Theorem
\ref{teoremaReduccion} this is equivalent to $\IA$ is a complete
intersection.

Assume that there exist $i,j: 1 \leq i < j \leq n$ such that $a_i,
a_j \in \mH$ and $m_i a_i = m_j a_j$, say $i = n-1$ and $j = n$, and
let us we prove that

\begin{center} {\rm CI-simplicial}$(\mA)$ = {\sc True} $\Longleftrightarrow$ $\IA$ is
a complete intersection. \end{center}

\noindent $(\Rightarrow)$ We prove that $\IA$ is a complete
intersection by induction on $n$. If $n = r+1$, we set $f =
x_{n-1}^{m_{n-1}} - x_n^{m_n}$. From the definition of $m_n$ we have
that $\gcd\{m_{n-1}, m_n\} = 1$; then the $1 \times n$ matrix whose
only row is $\widehat{f}$ is clearly dominating and  $\Delta_1 (A) =
\gcd\{m_{n-1}, m_n\} = 1$, thus by Theorem \ref{comprueba} it
follows that $\IA = (f)$ and $\IA$ is a complete intersection.

Assume that $n > r+1$ and that {\rm CI-simplicial}$(\mA)$ = {\sc
True}. First note that for every $V_l$ obtained during the execution
of {\rm CI-simplicial}$(\mA)$, all the vectors in $V_l$ are
proportional and $a_l = \gcd(V_l)$. Moreover, one observes that {\rm
CI-simplicial}$(V_l)$ = {\sc True}, hence by induction hypothesis
$I_{V_l}$ is a complete intersection whenever $V_l \subsetneq \mA$.

We consider be the set $\mB$ obtained after the {\bf while} loop,
then $\mB = \{a_{l_1},\ldots,a_{l_t}\}$ for some $1 \leq l_1 <
\cdots < l_t$, $t \geq 1$. We observe that $\mA = V_{l_1} \sqcup
\cdots \sqcup V_{l_t}$.

If $t = 1$, then $\mA = V_{l_1}$. We observe the last iteration of
the {\bf while} loop, where there exist $k_1, k_2: 1 \leq k_1 < k_2$
such that $m_{k_1} a_{k_1} = m_{k_2} a_{k_2} \in \N V_{k_1} \cap \N
V_{k_2}$ and $\mA = V_{l_1} = V_{k_1} \sqcup V_{k_2}$. We claim that
$\mA$ is a gluing of $V_{l_1}$ and $V_{l_2}$. Indeed, $m_{k_1}
a_{k_1} \in \N V_{k_1} \cap \N V_{k_2}$ and, by definition of
$m_{k_1}$ and $m_{k_2}$ it follows that $\Z V_{k_1} \cap \Z V_{k_2}
= \Z a_{k_1} \cap \Z a_{k_2} = \Z m_{k_1} a_{k_1}$. Hence $\IA$ is a
complete intersection by \cite[Theorem 1.4]{Rosalesgluing}.

If $t > 1$, then there exists $u \in \{1,\ldots,t\}$ such that
either:

(a) $a_{l_u} \in \Q (\mB \setminus \{a_{l_u}\})$ and setting
$B_{l_u} = {\rm min}\{b \in \Z^+ \, \vert \, b a_{l_u} \in \Z (\mB
\setminus \{a_{l_u}\})\}$ we have that $B_{l_u} a_{l_u} \in
\sum_{a_k \in V_{l_u}} \N a_k \cap \sum_{a_k \in \mA \setminus
V_{l_u}} \N a_k,$ or

(b) $a_{l_u} \notin \Q (\mB \setminus \{a_{l_u}\})$.

We denote $V := V_{l_u}$. In both cases we observe that
CI-simplicial$(\mA \setminus V) =$ CI-simplicial$(V) =$ {\sc True}.
Hence, by induction hypothesis, $I_{V}$ and $I_{\mA \setminus V}$
are complete intersections.

If (a) holds, then we have that $B_{l_u} a_{l_u} \in \N V \cap \N
(\mA \setminus V)$ and $\Z V \cap \Z (\mA \setminus V) = \Z a_{l_u}
\cap \Z (\mA \setminus V) = \Z B_{l_u} a_{l_u}$. Hence, $\mA$ is a
gluing of $V$ and $\mA \setminus V$ and $\IA$ is a complete
intersection by \cite[Theorem 1.4]{Rosalesgluing}.

If (b) holds, we denote $h_1 := {\rm ht}(I_{V})$, $n_1 := \vert V
\vert$, $h_2 := {\rm ht}(I_{\mA \setminus V})$ and $n_2 := n - n_1$
and take $\{f_1,\ldots,f_{h_1}\}$ and $\{g_1,\ldots,g_{h_2}\}$ two
minimal sets of generators of $I_{V}$ and $I_{\mA \setminus V}$
respectively. Hence, by Theorem \ref{comprueba} the matrices $A_1
\in \mathcal M_{h_1 \times n_1} (\Z)$ and $A_2 \in \mathcal M_{h_2
\times n_2}(\Z)$ whose rows are
$\widehat{f_1},\ldots,\widehat{f_{r_1}}$ and
$\widehat{g_1},\ldots,\widehat{g_{r_2}}$ respectively are both mixed
dominating and $\Delta_{h_1}(A_1) = \Delta_{h_2}(A_2) = 1$.
Moreover, ${\rm ht}(\IA) = h_1 + h_2$ and it follows that $\IA$ is a
complete intersection minimally generated by
$\{f_1,\ldots,f_{h_1},g_1,\ldots,g_{h_2}\}$. Indeed, the matrix $A
\in \mathcal M_{h \times n}$ whose rows are
$\widehat{f_1},\ldots,\widehat{f_{r_1}},\widehat{g_1},\ldots,\widehat{g_{r_2}}$
is clearly mixed dominating and $\Delta_h(A) = \Delta_{h_1}(A_1) \,
\Delta_{h_2}(A_2) = 1$.

\noindent $(\Leftarrow)$ Let us suppose that $\IA$ is a complete
intersection and $m_{n-1} a_{n-1} = m_n a_n$ and let us prove that
{\rm CI-simplicial}($\mA$) = {\sc True} by induction on $n$. If $n =
r+1$, then the result follows easily. Assume that $n > r + 1$. When
running {\rm CI-simplicial}$(\mA)$ we define: $G := \mA$, $V_k :=
\{a_k\}$ for all $k \in \{1,\ldots,n\}$ and $m_{k} := \min\{b \in
\Z^+\, \vert\, b  a_{k} \in \N\,(\mA \setminus \{a_{k}\})\}$ for all
$a_k \in {\rm Cone}(\mA \setminus \{a_k\})$.

 Since $m_{n-1}  a_{n-1} = m_n  a_n \in \N\, a_{n-1} \cap \N\, a_n$,
we denote $a_{n+1} := \gcd\{a_{n-1}, a_n\}$, $V_{n+1} :=
\{a_{n-1},a_n\}$ and $G := (G \setminus \{a_{n-1},a_n\}) \cup
\{a_{n+1}\}$; moreover, in case $a_{n+1} \in {\rm Cone}(G \setminus
\{a_{n+1}\})$, we define
\begin{center}$m_{n+1} := {\rm min}\{b \in \Z^+ \, \vert \, b
a_{n+1} \in \sum_{k = 1}^{n-2} \N\, a_k\}.$\end{center}

Consider now the set $\mA' := \{a_1,\ldots,a_{n-2},a_{n+1}\}$. By
Theorem \ref{ci-iguales} (a) we have that $I_{\mA^{\,\prime}}$ is a
complete intersection and by induction hypothesis it follows that
{\rm CI-simplicial}($\mA^{\,\prime}$) = {\sc True}. Now we run the
algorithm with
 $\mA^{\,\prime}$ as input and we get:
\begin{itemize}
\item $G^{\,\prime} := \mA^{\,\prime} = G$
\item $V_k^{\,\prime} := \{a_k\}$, for all $k \in \{1,\ldots,n-2,n+1\}$,
\item $m_k^{\,\prime} := {\rm min}\{b \in \Z^+ \, \vert \, b a_k \in \sum_{l \in
\{1,\ldots,n-2,n+1\} \atop l \neq k} \N\, a_l\}$ for all $k \in
\{1,\ldots,n-2,n+1\}$ such that $a_k \in {\rm Cone}(G^{\,\prime}
\setminus \{a_k\})$.
\end{itemize}
If $k \in \{1,\ldots,n-2\}$ then $V_k = V_k^{\,\prime}$ and if $a_k
\in {\rm Cone}(G^{\,\prime} \setminus \{a_k\})$, then by Theorem
\ref{ci-iguales} (b) we have that $m_k = m_k^{\,\prime}$; moreover
if $a_{n+1} \in {\rm Cone}(G^{\,\prime} \setminus \{a_{n+1}\})$,
then $m_{n+1}^{\,\prime} = m_{n+1}$ by definition.

Now, we continue with the execution of {\rm CI-simplicial}$(\mA)$
and {\rm CI-simplicial}$(\mA^{\,\prime})$ simultaneously. From now
on, we always have that $G^{\,\prime} = G$, hence $m_k' = m_k$ and
$V_k = V_k^{\,\prime}$ if $a_{n+1} \not\in V_k$ and $V_{k} =
(V_k^{\,\prime} \setminus \{a_{n+1}\}) \cup \{a_{n-1},a_n\}$ if
$a_{n+1} \in V_{k}^{\,\prime}$. Let us check that at any repetition
of the {\bf while} loop, whenever $m_k a_k = m_l a_l$ then $m_{k}
a_{k} \in \N\, V_{k}$.

Since {\rm CI-simplicial}$(\mA^{\,\prime})$ = {\sc True}, if $m_k
a_k = m_l a_l$, then $m_k a_k \in \N \,V_k^{\,\prime}$. If $a_{n+1}
\notin V_k^{\,\prime}$, then $V_k^{\,\prime} = V_k$ and we are done.
In case $a_{n+1} \in V_k^{\,\prime}$, then $V_l^{\,\prime} = V_l$,
 $m_k a_k = m_l a_l \in \N V_l$, and we deduce that
\begin{center}$ m_k = {\rm min}\{b \in \Z^+ \, \vert \, b  a_k \in \sum_{a_u
\in \mA' \setminus V_k'} \N\, a_u\}.$\end{center} We observe that we
are under the hypothesis of Lemma \ref{lematecnico} setting $i =
n-1$, $j = n$, $V' = V_k^{\,\prime}$ and $M = m_k$, thus $m_k a_k
\in \N V_k$. Proceeding analogously we get that $m_l a_l \in \N
V_l$.

Now we study the {\bf repeat} loop in the simultaneous execution. We
always have that $\mB^{\,\prime} = \mB$, then for all $a_k \in \mB$
we have that $B_k' = B_k$. Thus it only remains to prove that
whenever $B_k a_k \in \N V_k^{\,\prime} \cap \sum_{a_l \in (\mB
\setminus \{a_k\})} \N V_l^{\,\prime}$, then $B_k a_k \in \N V_k
\cap \sum_{a_l \in (\mB \setminus \{a_k\})} \N V_l$.

Note that if we denote $\mC := \cup_{a_l \in \mB} V_l$ we have that
$I_{\mC} \subset k[x_k \, \vert \, a_k \in \mC]$ is a complete
intersection by applying Lemmas \ref{redpertenece}, \ref{red2},
Proposition \ref{red} and Theorem \ref{ci-iguales} (a). Moreover, if
$a_{n-1}, a_n \notin \mC$, we have that $V_k^{\,\prime} = V_k$ for
all $a_k \in \mB$ and we can conclude that {\rm CI-simplicial}$(\mA)
=$ {\sc True}. In case $a_{n-1},\,a_n \in \mC$, we set $\tilde{m}_k
:= {\rm min}\{b \in \Z^+ \, \vert\, b a_k \in \sum_{a_r \in \mC
\setminus \{a_k\}} \N\, a_r\}$ for $k = n-1$ and $k = n$ and observe
that $\tilde{m}_{n-1}\, a_{n-1} = \tilde{m}_n \, a_n$ because
$m_{n-1} a_{n-1} = m_n a_n$.

Suppose that $B_k a_k \in \N\, V_k^{\,\prime}\, \cap \, \sum_{a_l
\in (\mB^{\,\prime} \setminus \{a_k\}) } \N\, V_l^{\,\prime}$, this
condition is equivalent to $B_k a_k \in \sum_{a_s \in
V_k^{\,\prime}} \N\, a_s \cap \sum_{a_r \in \mC^{\,\prime} \setminus
V_k^{\,\prime}} \N\, a_r$. Then
\begin{center} $B_k = {\rm min}\{b \in \Z^+ \, \vert\, b  a_k \in
\sum_{a_r \in \mC^{\,\prime} \setminus V_k^{\,\prime}} \N\, a_r\}.$
\end{center}

After applying Lemma \ref{lematecnico} to $I_{\mC}$, with $i = n-1$,
$j = n$, $V' = V_k^{\,\prime}$ and $M = B_k$, we have that $B_k a_k
\in \N V_k$ and that $B_k a_k \in \sum_{a_l \in \mC \setminus V_k}
\N \, a_l = \sum_{a_l \in (\mB \setminus \{a_k\})} \N V_l$. Thus, we
can conclude that CI-simplicial$(\mA) =$ {\sc True} and the proof is
complete. \QED \end{demo}

\medskip Let us illustrate how the algorithm CI-simplicial works with an
example.

\begin{Example}\label{ej2simplicial}Let us prove by means of {\rm Algorithm
CI-simplicial} that $\IA$ is a complete intersection, where $\mA :=
\{a_1,a_2,a_3,a_4,a_5,a_6,a_7,a_8\} \subset \N^3$ with $a_1 :=
(52,0,0)$, $a_2 := (0,52,0)$, $a_3 := (0,0,52)$, $a_4 :=
(20,30,100)$, $a_5 := (28,42,140)$, $a_6 := (30,45,150)$, $a_7 :=
(42,63,210)$ and $a_8 := (52,52,78)$.

We begin by setting $G := \mA$ and $V_i = \{a_i\}$ for all $i \in
\{1,\ldots,8\}$.
 For every $i \in \{4,5,6,7,8\}$ we observe that $a_i \in {\rm Cone}(\mA
\setminus \{a_i\})$, thus we compute $m_i$ and get that
\begin{center} $m_4 = 3$, $m_5 = 3$, $m_6 = 2$, $m_7 = 2$ and $m_8 =
2$. \end{center}

We observe that $3 a_4 = m_4 a_4 = m_6 a_6 = 2 a_6$. Then we set
$a_9 := \gcd\{a_4,a_6\} = (10,15,50)$, $G := (G \setminus
\{a_4,a_6\}) \cup \{a_9\} = \{a_1,a_2,a_3,a_5,a_7,a_8,a_9\}$ and
$V_9 := V_4 \cup V_6 = \{a_4,a_6\}$. Since $a_9 \in {\rm Cone}(G
\setminus \{a_9\})$, we also define $$m_9 := {\rm min}\{b \in \Z^+
\, \vert \, b a_9 \in \sum_{i \in \{1,2,3,5,7,8\}} \N a_i\} = 7.$$

We observe that $3 a_5 = m_5 a_5 = m_7 a_7 = 2 a_7$. Then we set
$a_{10} := \gcd\{a_5,a_7\} = (14,21,70)$, $G := (G \setminus
\{a_5,a_7\}) \cup \{a_{10}\} = \{a_1,a_2,a_3,a_8,a_9,a_{10}\}$ and
$V_{10} := V_5 \cup V_7 = \{a_5,a_7\}$. Since $a_{10} \in {\rm
Cone}(G \setminus \{a_{10}\})$, we also define $$m_{10} := {\rm
min}\{b \in \Z^+ \, \vert \, b a_{10} \in \sum_{i \in \{1,2,3,8,9\}}
\N a_i\} = 5.$$.

We observe that $7 a_9 = m_9 a_9 = m_{10} a_{10} = 5 a_{10}$ and we
check that $7 a_9 = 2 a_4 + a_6 = a_5 + a_7 \in \N V_9 \cap \N
V_{10}$. Then we set $a_{11} := \gcd\{a_9,a_{10}\} = (2,3,10)$, $G
:= (G \setminus \{a_{9},a_{10}\}) \cup \{a_{11}\} =
\{a_1,a_2,a_3,a_8,a_{11}\}$ and $V_{11} := V_{9} \cup V_{10} =
\{a_4,a_5,a_6,a_7\}$. Since $a_{11} \in {\rm Cone}(G \setminus
\{a_{11}\})$, we also define
$$m_{11} := {\rm min}\{b \in \Z^+ \, \vert \, b a_{11} \in \sum_{i
\in \{1,2,3,8\}} \N a_i\} = 52.$$

We observe that $m_i a_i \neq m_j a_j$ for every $a_i, a_j \in G$,
 then we take $\mB := G$.  If we denote $B_i := {\rm min}\{b\in \Z^+
\, \vert \, b a_i \in \Z (\mB \setminus \{a_i\})\}$ for all $i \in
\{1,2,3,8,11\}$, we get that $B_{11} = 52$. We observe that $B_{11}
a_{11} = 52 a_{11} = a_4 + 2 a_7 \in \N V_{11}$ and that $B_{11}
a_{11} = 52 a_{11} = a_2 + 7 a_3 + 2 a_8  \in \bigcup_{i \in
\{1,2,3,8\}} \N V_i$ and we define $\mB := \mB \setminus
\{a_{11}\}.$

Now we denote $B_i' := {\rm min}\{b\in \Z^+ \, \vert \, b a_i \in \Z
(G \setminus \{a_i\})\}$ for $i \in \{1,2,3,8\}$, and we get that
$B_8' = 2$. We observe that $B_8' a_8 = 2 a_8 \in \N V_{8}$ and that
$B_8' a_8 = 2 a_1 + 2 a_2 + 3 a_3 \in  \bigcup_{i \in \{1,2,3\}} \N
V_i$ and we define $\mB := \mB \setminus \{a_{8}\}.$

 We have that $\mB = \{a_1,a_2,a_3\}$, since
$a_1, a_2$ and $a_3$ are linearly independent we finally get that
$\mB = \emptyset$ and that {\rm IC-simplicial$(\mA)$} returns {\sc
True}. Consequently, $\IA$ is a complete intersection.
\end{Example}

\begin{Remark}\label{sistgen}The proofs of
{\rm Theorem \ref{algoritmosimplicial}} and \cite[Theorem
1.4]{Rosalesgluing} also show how to obtain a minimal set of
generators of $\IA$ when it is a complete intersection while
executing the {\rm Algorithm CI-simplicial}. Note that for obtaining
this minimal set of generators it is not necessary to perform any
additional calculations, but those of {\rm Algorithm CI-simplicial}.
More precisely, one can construct a minimal set of generators of
$\IA$ in the following way:
\begin{itemize}
\item[{\rm 1.}] While running the {\bf while} loop,
whenever there exist $a_i,a_j \in G$ such that $m_i a_i = m_j a_j$
we have that $m_i a_i \in \N\,V_i \cap \N\,V_j$. Hence, there exist
$\{\gamma_u \, \vert \, a_u \in V_i\}$, $\{\gamma_v\, \vert\, a_v
\in V_j \} \subset \N$ such that $m_i a_i = \sum_{a_u \in V_i}
\gamma_u a_u = \sum_{a_v \in V_j} \gamma_v a_v$ and we define the
binomial $g := \prod_{a_u \in V_i} x_u^{\gamma_u} - \prod_{a_v \in
V_j} x_v^{\gamma_v}$.

\item[{\rm 2.}] While running the {\bf repeat} loop, whenever $B_i a_i
\in \N\, V_i\, \cap\, \sum_{a_j \in \mB \setminus \{a_i\}} \N\,
V_j$, we have that there exist $\{\gamma_u\,\vert\,a_u \in V_i\}
\subset \N$ and
 $\bigcup_{a_j \in \mB \setminus \{a_i\}} \{\gamma_v\, \vert\, a_v \in V_j\} \subset \N$ such that $B_i a_i =
\sum_{a_u \in V_i} \gamma_u a_u = \sum_{a_j \in \mB \setminus
\{a_i\}}\left(\sum_{a_v \in V_j} \gamma_v a_v \right)$ and we define
the binomial $ g:= \prod_{a_u \in V_i} x_u^{\gamma_u} - \prod_{a_j
\in \mB \setminus \{a_i\}} (\prod_{a_v \in V_j} x_v^{\gamma_v})$.
\end{itemize}

\end{Remark}

\medskip Let us illustrate this procedure with {\rm Example
\ref{ej2simplicial}.

\begin{Example}We proved that $\IA$ is a complete intersection, where $\mA =
\{a_1,\ldots,a_8\}$ $\subset \N^3$ with $a_1 := (52,0,0)$, $a_2 :=
(0,52,0)$, $a_3 := (0,0,52)$, $a_4 := (20,30,100)$, $a_5 :=
(28,42,140)$, $a_6 := (30,45,150)$, $a_7 := (42,63,210)$ and $a_8 :=
(52,52,78)$.

\begin{itemize}
\item We computed $m_4,\, m_5,\, m_6,\, m_7$ and $m_8$ and observed that
$3 a_4 = m_4 a_4 = m_6 a_6 = 2 a_6$, then we define $g_1 := x_4^3 -
x_6^2$.

\item We defined $a_9 := \gcd\{a_4,a_6\} = (10,15,50)$,
computed $m_9$ and observed that $3 a_5 = m_5 a_5 = m_7 a_7 = 2
a_7$, then we define $g_2 := x_5^3 - x_7^2$.

\item We defined $a_{10} := \gcd\{a_5,a_7\} = (14,21,70)$, computed
$m_{10}$ and observed that  $7 a_9 = m_9 a_9 = m_{10} a_{10} = 5
a_{10}$ and checked that
 $m_9 a_9 = 2 a_4 + a_6 = a_5 + a_7 \in \N V_9 \cap \N V_{10}$,
then we define $g_3 := x_4^2 x_6 - x_5 x_7$.

\item Then we computed $B_{11}$ and checked that  $B_{11}
a_{11} = 52 a_{11} = a_4 + 2 a_7 = a_2 + 7 a_3 + 2 a_8 \in \N
\{a_4,a_5,a_6,a_7\} \cap \N \{a_1,a_2,a_3,a_8\}$, then we define
$g_4 := x_4 x_7^2 - x_2 x_3^7 x_8^2$.

\item Finally we computed $B_8'$ and checked that $B_{8}'
a_{8} = 2 a_8 = 2 a_1 + 2 a_2 + 3 a_3 \in \N \{a_8\} \cap \N
\{a_1,a_2,a_3\}$, then we define $g_5 := x_8^2 - x_1^2 x_2^2 x_3^3$.
\end{itemize}

The algorithm returned {\sc True}, then we can conclude that $\IA$
is a complete intersection minimally generated by
 $\{g_1,g_2,g_3,g_4,g_5\}$.
\end{Example}

\section{Ideal-theoretic complete intersection simplicial projective
 toric varieties and singularities.}\label{sec5}

This section is devoted to study the complete intersection property
in homogeneous simplicial toric ideals. The following result, which
is a direct consequence of Proposition \ref{simplicial}, shows how
Algorithm CI-simplicial can be simplified for homogeneous simplicial
toric ideals.

\begin{Corollary}\label{coralg}Let $\IA$ be a homogeneous simplicial toric
ideal. Then,
\begin{center} $\IA$ is a complete intersection $\Longleftrightarrow$  $\mA_{red} = \emptyset$.\end{center}
\end{Corollary}

This result yields an effective method for determining whether a
homogeneous simplicial toric ideal $\IA$ is a complete intersection.
Concretely, $\IA$ is a complete intersection if and only if the
algorithm in Table \ref{algRed} returns the empty set. Moreover, by
Remark \ref{generadores}, in case a homogeneous simplicial toric
ideal $\IA$ is a complete intersection, while checking that
$\mA_{red} = \emptyset$ one gets without any extra effort a minimal
set of generators of the toric ideal formed by binomials.

 For the rest of this section, $k$ denotes an algebraically
closed field and we aim at classifying those simplicial projective
toric varieties that are either smooth or have exactly one singular
point and are ideal-theoretic complete intersection. This
classification arises as a nontrivial consequence of Corollary
\ref{coralg}. We will prove the following two results, the first of
these is a direct consequence of Proposition \ref{2cartaslisas}.

\begin{Theorem}\label{smoothIC}Let $k$ be an algebraically closed field and $X \subset \P_k^{n-1}$ a smooth simplicial
 projective toric variety. Then, $X$ is an ideal-theoretic complete intersection if and only if
$n = 3$ and $X$ is the plane monomial curve defined parametrically
by $$x_1 = u_1^2 ,\, x_2 = u_2^2,\,  x_3 = u_1u_2.$$
 \end{Theorem}

\begin{Theorem}\label{1singIC}
Let $k$ be an algebraically closed field and $X \subset \P_k^{n-1}$
a simplicial projective toric variety with exactly one singular
point. Then, $X$ is an ideal-theoretic complete intersection if and
only if
\begin{itemize} \item either $X$ is the projective monomial curve in
$\P_k^{n-1}$ of degree $d \geq 3$ defined
 by
 $$x_1 = u_1^d,\,  x_2 = u_2^d,\, x_3 = u_1^{d-1} u_2,\, x_4 = u_1^{d-d_4} u_2^{d_4} \ldots,\,  x_{n} =
 u_1^{d-d_n} u_2^{d_n},$$ where $1 < d_4 < \cdots < d_n < d$ and $d_4 \mid
 d_5 \mid \cdots \mid d_n \mid d$,
 \item or $X$ is the projective monomial surface in $\P_k^3$
 defined by $$x_1 = u_1^2,\, x_2 = u_2^2,\, x_3 = u_3^2,\, x_4 = u_1u_2.$$
 \end{itemize}
 \end{Theorem}

To obtain these results, we study the affine pieces of a simplicial
projective toric variety. Let $X$ be a simplicial projective toric
variety, then there exists a set $\mA = \{d e_1, \ldots, d e_m,
a_{m+1}, \ldots, a_n\} \subset \N^m$ such that $X = V(\IA)$, where
$\{e_1,\ldots,e_m\}$ is the canonical basis of $\Z^m$ and $d =
\sum_{j = 1}^m a_{ij} \in \Z^+$ for all $i \in \{m+1,\ldots,n\}$
(see, e.g., \cite[Section 2]{HerzogHibi}).

Consider the affine pieces $\{X \cap \mathcal U_i\}_{i = 1}^n$ of
$X$, where $\mathcal U_i := \P_k^{n-1} \setminus V(x_i)$ for all $i
\in \{1,\ldots,n\}$. Since $X$ is simplicial, it suffices to
consider the $m$ first affine pieces for covering $X$. This is, $X =
\bigcup_{i = 1}^m\ (X \cap \mathcal U_i).$ Indeed, if $p = (p_1 :
\cdots : p_n) \in X$ and $p \notin U_i$ for all $i \in
\{1,\ldots,m\}$, then $p_1 = \cdots = p_m = 0$. Now for all $j \in
\{m+1,\ldots,n\}$, we consider the binomial $f_j := x_j^d - \prod_{k
= 1}^m x_k^{a_{j k}} \in \IA$. Since $f_j(p) = 0$, we deduce that
$p_j = 0$ for all $j \in \{m+1,\ldots,n\}$, which is not possible.

Let us see that these affine pieces are homeomorphic to affine
simplicial toric varieties. Recall that for all $i \in
\{1,\ldots,n\}$, $U_i \simeq \A_k^{n-1}$ via
$$(b_1:\ldots:b_n) \mapsto (b_1/b_i, \ldots,\, b_{i-1}/b_i,\,
b_{i+1}/b_i,\ldots,\, b_n/b_i).$$ For all $i \in \{1,\ldots,m\}$ and
for all $j \in \{m+1,\ldots,n\}$ we denote
 $$a_j^{(i)} := (a_{j 1},\ldots,a_{j
i-1},a_{j  i+1},\ldots,a_{j m}) \in \N^{m-1}$$ and $\mathcal A^{(i)}
:= \{d \bar{e}_1,\ldots, d \bar{e}_{m-1},
a_{m+1}^{(i)},\ldots,a_n^{(i)}\} \subset \N^{m-1}$ where
$\{\bar{e}_1,\ldots,\bar{e}_{m-1}\}$ is the canonical basis of
$\Z^{m-1}$, then the affine piece $X \cap \mathcal U_i$ is
homeomorphic to the affine simplicial toric variety  $Y_i :=
V(I_{\mA^{(i)}}).$

Thus, $X$ is smooth if and only if $Y_i$ is smooth for all $i \in
\{1,\ldots,m\}$. The following classical result characterizes when
an affine toric variety is smooth.

\begin{Theorem}\citep{CLS, Fulton, KKMS}\label{smooth}
Let $k$ be an algebraically closed field. Then the following
conditions are equivalent:
\begin{itemize}
\item $V(\IA)$ is smooth.
\item $0 \in \A_k^n$ is a regular point of $V(\IA)$.
\item The semigroup $\sum_{i = 1}^n \N a_i$ admits a set of generators with ${\rm dim}(\Q \mA)$ elements.
\end{itemize}
\end{Theorem}

In particular we have that every smooth affine toric variety is
simplicial and have the following corollary, whose proof is easy.

\begin{Corollary}\label{simplicialafinlisa} Let $k$ be an algebraically closed field and $X \subset \A_k^n$
an affine toric variety. Then, $X$ is smooth $\Longleftrightarrow$
$X = V(\IA)$ where $\mA = \{d e_1,\ldots,d e_m, a_{m+1},\ldots,a_n\}
\subset \N^m$ and for all $j \in \{1,\ldots,m\}$, if $\lambda_j :=
{\rm min}\{k \in \Z^+ \, \vert \, k e_j \in \mA\}$, then  $\lambda_j
\mid d$ and $\lambda_j \mid a_{ij}$ for all $i \in
\{m+1,\ldots,n\}.$
\end{Corollary}

We are making use of Corollary \ref{simplicialafinlisa} to prove the
next proposition, which is a consequence of Corollary \ref{coralg}.
From Proposition \ref{2cartaslisas} the proof of Theorem
\ref{smoothIC} follows at once, moreover it is useful to prove
Theorem \ref{1singIC}.

\begin{Proposition}\label{2cartaslisas}Let $k$ be an algebraically closed field and
 $X = V(\IA) \subset \P_k^{n-1}$ $(n \geq 3)$ a simplicial projective toric
variety with $\mA = \{de_1,\ldots,d e_m, a_{m+1},\ldots,a_n\}
\subset \N^m$ and $\sum_{j = 1}^n a_{ij} = d$ for all $i \in
\{m+1,\ldots,n\}$. Suppose that there exist $r,s: 1 \leq r < s \leq
m$ such that every point in $X \setminus V(x_r,x_s)$ is a regular
point. Then, $X$ is an ideal-theoretic complete intersection if and
only if $n = m+1$ and $X$ is given parametrically by $$x_1 =
u_1^2,\, \ldots,\, x_m = u_m^2, \, x_{m+1} = u_i u_j. $$
\end{Proposition}
\begin{demo}Suppose that $\gcd\{d, a_{ij}\, \vert \, m+1\leq i \leq n,\, 1
\leq j \leq m\} = 1$ and that $r = 1$ and $s = 2$. Then, $Y_k =
V(I_{\mA^{(k)}})$ is smooth for $k \in \{1,2\}$. We first prove that
$e_1 + (d-1)e_2, (d-1) e_1 + e_2 \in \mA$. We denote
$$\lambda_{1j} := {\rm min} \{k \in \Z^+\, \vert\, (d-k) e_1 + k e_j \in \mA\}
{\rm \ for\  all\ } j \in \{2,\ldots,m\} {\rm, \ and}$$
$$\lambda_{2j} := {\rm min} \{k \in \Z^+\, \vert\, (d-k) e_2 + k e_j \in \mA\}
{\rm \ for\  all\ } j \in \{1,3,\ldots,m\}.$$

By Corollary \ref{simplicialafinlisa} it follows that for all $k \in
\{1,2\}$
\begin{center} $\lambda_{k j} \mid d \ {\rm and}
\ \lambda_{k j}\, |\, a_{i j}$ for all $i \in \{m+1,\ldots,n\}$ and
for all $j \in \{1,\ldots,m\} \setminus \{k\}.$
\end{center}

We claim that $\lambda_{1 2} = \lambda_{2 1} = 1$. Indeed,
$\lambda_{1 2} \mid d$ and  $\lambda_{2 1} e_1 + (d - \lambda_{2
1})e_2 \in \mA$, then $\lambda_{1 2} \mid d - \lambda_{2 1}$, which
implies that $\lambda_{1 2} \mid \lambda_{2 1}$ and by a similar
argument we get that $\lambda = \lambda_{2 1}$. Moreover, for all $j
\geq 3$, $(d - \lambda_{1 j}) e_1 + \lambda_{1 j} e_j \in \mA$,
which implies that $\lambda_{1 2} \mid \lambda_{1 j}$ and
consequently $\lambda_{1 2} \mid \gcd\{d, a_{ij}\, \vert \, m+1\leq
i \leq n,\, 1 \leq j \leq m\} = 1$. Hence, $\lambda_{1 2} =
\lambda_{2 1} = 1$ and $(d-1)e_1 + e_2,\, e_1 + (d-1) e_2 \in \mA$.

If $d > 2$, we may assume that $a_{n-1} = e_1 + (d-1)e_2$ and $a_n =
(d-1)e_1 + e_2$. Equality $a_{n-1} + a_n = d e_1 + d e_2$ implies
that $d e_1, d e_2, a_{n-1}, a_n \in \mA_{red}$ and by Corollary
\ref{coralg} we get hat
 $\IA$ is not a complete intersection.

 If $d = 2$, we have proved that $e_1+e_2 \in \mA$. If $n =
m+1$ we have that $\mA = \{ 2 e_1,\ldots, 2 e_m, e_1 + e_2\}$. In
this case $B_n = 2$, $2 a_n = 2 e_1 + 2 e_2$ and $\{2 e_1,\ldots,2
e_m\}$ are linearly independent, thus $\mA_{red} = \emptyset$ and
$\IA$ is a complete intersection by Corollary \ref{coralg}. If $n >
m+1$, then there exists $e_i + e_j \in \mA$ with $j \geq 3$, $i \neq
j$, which implies that $\lambda_{1 j} = \lambda_{2 j} = 1$. Hence
$e_1 + e_j$ and $e_2 + e_j$ belong to $\mA$. Since $e_1 + e_2 , 2
e_1, 2 e_2, 2 e_j \in \mA$, we have the following equality involving
six elements of $\mA$
$$(e_1 + e_j) + (e_2 + e_j) + (e_1 + e_2) = 2 e_1 + 2
e_2 + 2 e_j.$$ From this equality we derive that $\mA_{red}$ has at
least six elements and cannot be empty. By Corollary \ref{coralg},
we get that $\IA$ is not a complete intersection and the result
follows. \QED \end{demo}

\vspace{1cm}

\noindent {\bf Proof of Theorem \ref{1singIC}.} Take $\mA = \{d
e_1,\ldots,d e_m, a_{m+1},\ldots,a_n\} \subset \N^m$ with $d =
\sum_{j = 1}^m a_{ij}$ for all $i \in \{m+1,\ldots,n\}$ such that $X
= V(\IA)$. We claim that there exists an only affine piece $Y_i$
with $i \in \{1,\ldots,m\}$ which is not smooth. Indeed, if there
exist $i,j: 1 \leq i < j \leq m$ such that $Y_i, Y_j$ are not
smooth, say $Y_1$ and $Y_2$, then by Theorem \ref{smooth}, $0 \in
Y_1 \subset \A_k^{n-1}$ is a singular point of $Y_1$ which
corresponds to the singular point  $(1:0:\cdots:0) \in X$ and $0 \in
Y_2 \subset \A_k^{n-1}$ is a singular point of  $Y_2$ which
corresponds to the singular point  $(0:1:0:\cdots:0) \in X$, a
contradiction. Then we can assume that the affine pieces
$Y_1,\ldots,Y_{m-1}$ are smooth. If $m \geq 3$, by Proposition
\ref{2cartaslisas}, we directly get that $m = 3$, $n = 4$ and $X$ is
given parametrically by
$$x_1 = u_1^2,\, x_2 = u_2^2,\, x_3 =  u_3^2,\, x_4 = u_1 u_2.$$

Thus, it only remains to consider the case of projective monomial
curves with exactly a singular point. In this setting we have that
 $a_3 = (d - d_3) e_1 + d_3 e_2, \ldots, a_n
= (d - d_n) e_1 + d_n e_2$ and we can assume that $d_3 < \cdots <
d_n < d$ and that  $\gcd\{d, d_3, \ldots, d_n\} = 1$. We have that
$Y_1 = V(I_{\mA^{(1)}})$ where $\mA^{(1)} = \{d,d_3,\ldots,d_n\}$.
By Corollary \ref{simplicialafinlisa}, we have that  $Y_1$ is smooth
if and only if $d_3 = \min (\mA^{(1)})$, divides to every element of
$\mA^{(1)}$. Since $\gcd\{d,d_3,\ldots,d_n\} = 1$, then $d_3 = 1$.

 Now we prove by induction on $n \geq 3$ that $\IA$
with
$$\mA = \{d e_1,\, d e_2,\, (d - 1) e_1 +  e_2, (d - d_4) e_1 + d_4
e_2, \ldots,\, (d - d_n) e_1 + d_n e_2\}$$  is a complete
intersection if and only if $d_4 \mid \cdots \mid d_n \mid d$.

 If $n = 3$ we have that $\mA = \{de_1, de_2, (d-1) e_1 +
e_2\}$. In this setting we have that $B_3 = d$ and $d a_3 = (d-1) d
e_1 + d e_2$, being $B_3 = {\rm min}\{b \in \Z^+ \, \vert \, b a_3
\in \Z d e_1 + \Z d e_2\}$. Since $d e_1, d e_2$ are linearly
independent we obtain that $\mA_{red} = \emptyset$ and by Corollary
\ref{coralg}, $\IA$ is a complete intersection.

 Assume now that $n \geq 4$, we observe that $a_i = d_i
a_3 + (1 - d_i) d e_1$ for all $i \in \{4,\ldots,n\}$, hence $B_i =
1$. Moreover, it is easy to check that $b a_3 \in \sum_{i \in
\{1,\ldots,n\} \atop i \neq 3} \Z a_i$  if and only if $b \in
\sum_{i = 4}^n \Z d_i + \Z d$, from where we get that $B_3 =
\gcd\{d_4,\ldots,d_n,d\}$. Moreover
\begin{center}$B_3 a_3 \in \sum_{j \in \{4,\ldots,n\}} \N
a_j + \N d e_1 + \N d e_2 \Longleftrightarrow d_4 =
\gcd\{d_4,\ldots,d_n,d\},$\end{center} indeed if $B_3 a_3 \in
\sum_{j \in \{4,\ldots,n\}} \N a_j + \N d e_1 + \N d e_2$, then we
have that $\gcd\{d_4,\ldots,d_n,d\} \in \sum_{j = 4}^{m} \N d_j + \N
d$, but $\gcd\{d_4,\ldots,d_n,d\} \leq d_4 < \cdots < d_n < d$, so
this can only happen if $d_4 = \gcd\{d_4,\ldots,d_n,d\}$. Note that
in this case we have that $B_3 a_3 = a_4 + (d_4 - 1) d e_1$.

 Hence, if $d_4 \not= \gcd\{d_4,\ldots,d_n,d\}$,
$\mA_{red} \not= \emptyset$ and $\IA$ is not a complete intersection
by Corollary \ref{coralg}. Otherwise, i.e., if $d_4 =
\gcd\{d_4,\ldots,d_n,d\}$, by Proposition \ref{red} and Lemma
\ref{redpertenece}, we have that $\IA$ is a complete intersection if
and only if $I_{\mA \setminus \{a_3\}}$ so is. In this case,
denoting $d^{\,\prime} := d / d_4 \in \Z^+$, $d_i^{\,\prime} := d_i
/ d_4 \in \Z^+$ for all $i \in \{4,\ldots,n\}$ and $\mA^{\,\prime}
:= \{d^{\,\prime} e_1, d^{\,\prime} e_2, (d^{\,\prime} - 1) e_1 +
e_2, (d^{\,\prime} - d_5^{\,\prime}) e_1 + d_5^{\,\prime} e_2,
\ldots,\, (d^{\,\prime} - d_n^{\,\prime}) e_1 + d_n^{\,\prime} e_2
\}$, it is evident that $I_{\mA \setminus \{a_3\}} =
I_{\mA^{\,\prime}}$. By induction hypothesis $I_{\mA^{\,\prime}}$ is
a complete intersection if and only if $d_5^{\,\prime} \mid \cdots
\mid d_n^{\,\prime} \mid d^{\,\prime}$, which implies that $\IA$ is
a complete intersection $\Longleftrightarrow\, d_4 \mid d_5 \mid
\cdots \mid d_n \mid d$. \QED

Note that in the proofs of Theorems \ref{smoothIC}, \ref{1singIC}
and Proposition \ref{2cartaslisas} we have applied the algorithm in
Table \ref{algRed} to check whether $\mA_{red} = \emptyset$. Hence,
following Remark \ref{generadores}, we get the defining equations of
the ideal-theoretic complete intersection simplicial projective
toric varieties that are smooth or have one singular point. More
precisely, we get the following results.

\begin{Corollary}Let $k$ be an algebraically closed field and $X \subset \P_k^{n-1}$ a smooth simplicial
 projective toric variety. Then, $X$ is an ideal-theoretic complete intersection if and only
 if $X$ is the curve in $\P_k^2$ with equation $x_3^2 - x_1 x_2 = 0$.
\end{Corollary}

\begin{Corollary}
Let $k$ be an algebraically closed field and $X \subset \P_k^{n-1}$
a simplicial projective toric variety with exactly one singular
point. Then, $X$ is an ideal-theoretic complete intersection if and
only if
\begin{itemize} \item either $X$ is the monomial curve in
$\P_k^{n-1}$, with defining equations

$$\left\{ \begin{array}{ccc} x_3^{b_3} - x_1^{b_3 - 1} x_4 & = & 0
\\ x_4^{b_4} - x_1^{b_4 - 1}  x_5& = &0 \\ \vdots \\
x_{n-1}^{b_{n-1}} - x_1^{b_{n-1} - 1} x_n& =& 0\\
x_n^{b_n} - x_1^{b_n - 1} x_2 &= &0\end{array} \right.$$
 where either $n \geq 4$ and $b_3,\ldots,b_n \geq 2$,  or $n = 3$ and
$b_3 \geq 3$,

 \item or $X$ is the surface in $\P_k^3$ of degree $2$ with equation $x_4^2 - x_1
 x_2 = 0$.
 \end{itemize}
 \end{Corollary}

\section{Computational aspects}\label{sec6}

In this section we explain how we have implemented the algorithms
obtained in Section 4 and 5 for checking whether a simplicial toric
ideal or a homogeneous simplicial toric ideal is a complete
intersection. We have implemented these algorithms in C++ and in
{\sc Singular}. The implementation in {\sc Singular} gave rise to
the distributed library {\tt cisimplicial.lib}, which is included in
the software since its version 3-1-4.

Given $\IA$ a homogeneous simplicial toric ideal, according to
Corollary \ref{coralg}, to check whether $\IA$ is a complete
intersection one can verify if $\mA_{red} = \emptyset$. For this
purpose one has to design procedures to solve the following problems
\begin{itemize}
\item to check whether an element $b \in \N^m$ belongs to a subsemigroup of $\N^m$

\item to compute $B_i := {\rm min}\{b \in \Z^+\, \vert \ b\, a_i \in
 \sum_{j \in \{1,\ldots,n\} \atop j \neq i} \Z\, a_j \}$ for all $i \in
 \{1,\ldots,n\}$.
\end{itemize}

For solving the first problem we have implemented a full enumeration
procedure. An algorithm based on a Graph Theory approach can be
found in \cite{CD}, for yet another methods we refer the reader to
\cite{PV} and the references there.

\begin{Example} In {\tt cisimplicial.lib} we have implemented the function
{\rm \textbf{belongSemigroup}}, which performs an enumeration to
check whether an element $b \in \N^m$ belongs to a subsemigroup of
$\N^m$. This function receives as input an $m \times n$ integral
matrix $A$ and a vector $b \in \N^m$ and checks whether the system
of equations $A\, x = b$ has a solution $x \in \N^n$, it returns a
solution whenever such it exists or $0$ if it does not exist. The
following example shows how to use {\tt cisimplicial.lib} to check
that $b := (22,12,10)$ belongs to the semigroup spanned by $\mA$,
meanwhile $c := (12,4,1)$ does not, where $\mA = \{(10,2,5),\,
(3,1,0),\, (2,1,1),\, (1,3,2)\}$.

\medskip
\parskip 2pt
{\rm \texttt{$>$ intmat A[3][4] $=$ 10,\ 3,\ 2,\ 1,}

\texttt{\ \ \ \ \ \ \ \ \ \ \ \ \ \ \ \ \ \ \ \ \ 2,\ 1,\ 1,\ 3,}

\texttt{\ \ \ \ \ \ \ \ \ \ \ \ \ \ \ \ \ \ \ \ \ 5,\ 0,\ 1,\ 2;}

\texttt{$>$ intvec b $=$ 23,12,10;}

\texttt{$>$ intvec c $=$ 12,4,1;}

\texttt{$>$ belongSemigroup(b,A);}

\texttt{ 1,3,1,2 \hspace{1cm} //  A * (1 3 1 2)$^{\texttt t}$ = b}

\texttt{$>$ belongSemigroup(c,A);}

\texttt{ 0\ \ \ \ \ \ \ \hspace{1cm} // A * x$^{\texttt t}$ = c has
no solution x $\in \N^n$}}
\end{Example}

\parskip 8pt
To solve the second problem we use the following result:
\begin{Lemma}\label{lemabi}Let $\mA = \{a_1,\ldots,a_n\}$ be a set of vectors in $\N^m$ and suppose that
$a_i \in \sum_{j \in \{1,\ldots,n\} \atop j \neq i} \Q a_j$ for some
$i \in \{1,\ldots,n\}$. Then
$$
B_i = \frac{{\rm Card} (T(\Z^m / \sum_{j \in \{1,\ldots,n\} \atop j
\neq i} \Z a_j)) }{ {\rm Card} (T(\Z^m / \sum_{j \in \{1,\ldots,n\}}
\Z a_j))},
$$
where $T(-)$ denotes de torsion subgroup of an abelian group.
\end{Lemma}
\begin{demo}Consider the group $G := \Z^m / \sum_{j \in \{1,\ldots,n\} \atop j
\neq i} \Z a_j$,  then $a_i \in T(G)$ and $B_i$ is the order of
$a_i$ in $G$. Thus $B_i = {\rm Card} (T(G)) / {\rm Card} (T(\Z^m /
\sum_{j \in \{1,\ldots,n\}} \Z a_j))$. \QED \end{demo}

 This result reduces the problem of computing $B_i$ to that
of computing the order of the torsion subgroup of two finitely
generated abelian groups, which we compute in polynomial time by
means of the Hermite Normal Form of a matrix, see \cite[Theorem
2.4.3 and Algorithm 2.4.5]{Cohen}.

\begin{Example} In {\tt cisimplicial.lib} we compute $B_i$ by means of the
function {\rm \textbf{cardGroup}}. This function computes the order
of a finite abelian group, it receives as input an $m \times n$
integral matrix $C$ and returns the order of the group $\Z^n / G$,
or $-1$ if it is infinite, where $G$ denotes the group spanned by
the columns of $C$. The following example shows how to compute with
{\tt cisimplicial.lib} the order of $\Z^3 / G$, where $G$ is the
group spanned by $\{(24,0,0),\, (0,24,0),\, (0,0,24),\, (8,10,5),\,
(3,6,9)\}$.

\bigskip

\parskip 2pt
{\rm

\texttt{$>$ intmat C[3][5] $=$ 24,\ \ 0,\ \ 0,\ \ 8,\ \ 3,}

\texttt{\ \ \ \ \ \ \ \ \ \ \ \ \ \ \ \ \ \ \ \ \ 0,\ 24,\ \ 0,\
10,\ \ 6,}

\texttt{\ \ \ \ \ \ \ \ \ \ \ \ \ \ \ \ \ \ \ \ \ 0,\ \ 0,\ 24,\ \
5,\ \ 9;}

\texttt{$>$ cardGroup(C);}

\texttt{\ 72\ \ \ \ \ \ \ \hspace{1cm} // The order of the group}
$\Z^3 / G$ \texttt{is} $72$}
\end{Example}
\parskip 8pt

 To implement the {\rm Algorithm CI-simplicial}, besides the
problems of computing $B_i$ and checking whether a vector belongs to
a subsemigroup of $\N^m$, a key point is the computation of $m_i$
for every $a_i$ such that $a_i \in {\rm Cone}(\mA \setminus
\{a_i\})$. If one wants to compute explicitly these values one could
generalize the method introduced in \cite[Section 4.1.2]{BGS} based
on a Graph Theory representation of the problem. However, we have
perfomed an implementation in which we do not aim at finding the
optimum value $m_i$ but only at checking whether there exist
$a_i,\,a_j$ such that $m_i a_i  = m_j a_j$. For this purpose, for
every $a_i \in {\rm Cone}(\mA \setminus \{a_i\})$ we define the set
$\mathcal L_i := \{a_j \in \mA \,\vert\, \exists \lambda \in \Q$
such that $\lambda a_i = a_j\}$, take $k(i) \in \{1,\ldots,m\}$ such
that $a_{i\,k(i)} \neq 0$ and set
\begin{center}$\overline{m}_i := {\rm min}\{b \in \Z^+\,\vert\, b\, a_{i\,k(i)}
\in \sum_{a_j \in \mathcal L_i \atop j \neq i} \N\,
a_{j\,k(i)}\}.$\end{center} Note that $\overline{m}_i$ can be
computed following the method proposed in \cite[Section 4.1.2]{BGS}.
Moreover, the following properties hold:
\begin{itemize} \item $\overline{m}_i \geq m_i$
because $\overline{m}_i a_i \in \sum_{a_j \in \mathcal L_i \atop j
\neq i} \N a_j$
\item  if $m_i a_i = m_j a_j$, then $a_j \in \mathcal L_i$ and $m_i =
\overline{m}_i$.
\item $m_i = \overline{m}_i$ if and only if the only solution $(x_1,\ldots,x_n) \in
\N^n$ to the system of equations $x_1 a_1 + \cdots + x_n a_n =
(\overline{m}_i - 1) a_i$ is the trivial one, i.e.,  $x_i =
\overline{m}_i - 1$ and $x_j = 0$ for every $j \in \{1,\ldots,n\}
\setminus \{i\}$. In particular, if $\overline{m}_i = 1$, then $m_i
= 1$.
 \end{itemize}

Thus, one can avoid the exact computation of $m_i$ by computing
$\overline{m}_i$ and by checking whether a system of diophantine
equations has more than one nonnegative integral solution. Moreover,
this problem is equivalent to checking whether the system of
equations
\begin{center} $\left. \begin{array}{c} x_1 a_1 + \cdots + x_n a_n = (\overline{m}_i - 1) a_i \\
x_i + x_{n+1} = \overline{m}_i - 2   \end{array} \right\} $
\end{center}
has a solution $(x_1,\ldots,x_n,x_{n+1}) \in \N^{n+1}$, which we
solve by enumeration.

We have produced an implementation of Algorithm CI-simplicial
following the techniques we have already described. Computational
experiments show that our implementation of CI-simplicial is able to
solve large size instances. For example, we have produced examples
of simplicial toric ideals $\IA$ with $\mA = \{a_1,\ldots,a_{27}\}
\subset \N^8$ and $0 \leq a_{ij} \leq 4000$ for all $1 \leq i \leq
27$, $1 \leq j \leq 8$ and we have determined in less than a second
on a personal computer with Intel Pentium IV 3Ghz whether $\IA$ is a
complete intersection.

\begin{Example}The main function in {\tt cisimplicial.lib} is called {\rm \textbf{isCI}}.
It receives an $m \times n$ integral matrix $A$ whose columns
correspond to the vectors in the set $\mA = \{d_1 e_1,\ldots,d_m e_m
, a_{m+1},\ldots,a_n\} \subset \N^m$ and returns $0$ if the
simplicial toric ideal $\IA$ is not a complete intersection or $1$
otherwise. Whenever $\IA$ is a complete intersection, it also
returns a minimal set of binomials generating the ideal. The
following example shows how to use {\tt cisimplicial.lib} to check
that $\IA$ is not a complete intersection and $I_{\mB}$ is a
complete intersection, where:
\begin{itemize} \item $\mA := \{(12,0,0),\, (0,10,0),\, (0,0,8),\, (1,3,3),\,
(2,2,3)\}$, and \item $\mB :=
\{(52,0,0),\,(0,52,0),\,(0,0,52),\,(20,30,100),\,(28,42,140),\,(30,45,150),\\
\,(42,63,210),\, (32,32,48),\,(36,36,54),\,(40,40,6012)\}$.
\end{itemize}

 {\rm

  \medskip

\parskip2pt

\texttt{$>$ intmat A[3][5] $=$ 12,\ \ 0,\ \ 0,\ \ 1,\ \ 2,}

\texttt{\ \ \ \ \ \ \ \ \ \ \ \ \ \ \ \ \ \ \ \ \ 0,\ 10,\ \ 0,\ \
3,\ \ 2,}

\texttt{\ \ \ \ \ \ \ \ \ \ \ \ \ \ \ \ \ \ \ \ \ 0,\ \ 0,\ \ 8,\ \
3,\ \ 3;}

\texttt{$>$ isCI(A);}

\texttt{\ 0 \hspace{3cm} // It is not a complete intersection}

\texttt{$>$ intmat B[3][10] $=$}

\texttt{\ 52,\ \ \ 0,\ \ \  0,\ \  20,\ \ 28,\ \ 30,\ \ 42,\ \ 32,\
\ 36,\ \ 40,}

\texttt{\ \ 0,\ \ 52,\ \ \ 0,\ \ 30,\ \ 42,\ \ 45,\ \ 63,\ \ 32,\ \
36,\ \ 40,}

\texttt{\ \ 0,\ \ \ 0,\ \ 52,\  100,\  140,\  150,\  210,\ \ 48,\ \
54, 6012;}

\texttt{$>$ isCI(B);}

\texttt{ 1 \hspace{3cm} // It is a complete intersection}

\texttt{ // Generators of the toric ideal}

\texttt{ toric[1]=x(4)\^{ }3-x(6)\^{ }2}

\texttt{ toric[2]=x(5)\^{ }3-x(7)\^{ }2}

\texttt{ toric[3]=x(4)\^{ }2*x(6)-x(5)*x(7)}

\texttt{ toric[4]=x(8)\^{ }5-x(10)\^{ }4}

\texttt{ toric[5]=-x(9)\^{ }2+x(8)*x(10)}

\texttt{ toric[6]=-x(1)\^{ }2*x(2)\^{ }3*x(3)\^{ }10+x(4)*x(7)\^{
}2}

\texttt{ toric[7]=-x(1)\^{ }2*x(2)\^{ }2*x(3)\^{ }3+x(8)\^{
}2*x(10)}}

\end{Example}

\parskip 8pt

%\bibliographystyle{elsart-harv}
%\bibliography{mybibfile}

% Include the ".bib" file (generated by bibtex) right here.

\end{document}